\documentstyle[12pt,epic,eepic]{article}

\newfont{\goth}{eufm10 scaled \magstep1}
\newtheorem{example}{Example}[section]
\newtheorem{remark}[example]{Remark}
\newtheorem{theorem}[example]{Theorem}
\newtheorem{corollary}[example]{Corollary}

\newtheorem{proposition}[example]{Proposition}

\newtheorem{lemma}[example]{Lemma}
\newtheorem{conjecture}[example]{Conjecture}

\def\Proof{\medskip\noindent {\it Proof --- \ }}
\def\cqfd{\hfill $\Box$ \bigskip}
\def\adots{\mathinner{\mkern2mu\raise1pt\hbox{.}
 \mkern3mu\raise4pt\hbox{.}\mkern1mu\raise7pt\hbox{.}}}
\def\<{\langle}
\def\>{\rangle}

\def\cf{{\it cf. }}
\def\ie{{\it i.e. }}

\def\col{{\rm col\,}}
\def\yam{{\rm Yam}}

\def\S{{\goth S}}
\def\sym{{\sl Sym}}
\newfont{\bb}{cmbx10}
\def\N{{\bf N}}
\def\Z{{\bf Z}}

\def\Q{{\bf Q}}

\def\Q{{\bf Q}}

\def\mod{{\rm \ mod\ }}
\def\C{{\bf C}}

\def\P{{\cal P}}

\def\tab{{\rm Tab\, }}

\def\h{\hbox{\goth h}}

\def\F{{\cal F}}

\def\H{{\cal H}}
\def\dim{{\rm dim\,}}

\def\F{{\cal F}}

\def\tab{{\rm Tab\,}}
\def\t{{\bf t}}
\def\tQ{{\tilde Q'}}
\def\tK{{\tilde K}}
\def\tG{{\tilde G}}
\def\tH{{\tilde H}}
\def\GF{{\bf F}}
\def\cc{{\tilde c}}
\def\cs{{\tilde s}}
\def\G{{\cal G}}

\def\Sl{\mbox{\goth sl}}
\def\gl{\mbox{\goth gl}}
\def\slchap{\widehat{\Sl}}
\def\sym{{\sl Sym}}
\def\H{{\cal H}}
\def\T{{\cal S}}
\def\h{\mbox{\goth h}}
\def\F{{\cal F}}

\def\wed{\wedge_q}
\def\bigwed{\bigwedge_q}

\def\tA{\tilde A}

\title{RIBBON TABLEAUX,\\ HALL-LITTLEWOOD FUNCTIONS,\\
QUANTUM AFFINE ALGEBRAS
\\ AND UNIPOTENT VARIETIES \thanks{Partially supported
by PRC Math-Info and EEC grant n$^0$ ERBCHRXCT930400} }

\author{Alain {\sc Lascoux}\thanks{Institut Blaise Pascal, L.I.T.P.,
Universit\'e
Paris 7, 2 place Jussieu, 75251 Paris cedex 05, France},
\rm Bernard {\sc  Leclerc}\thanks{Universit\'e de Caen, D\'epartement
de Math\'ematiques, Esplanade de la Paix, BP 5186, 14032 Caen cedex, France}
\rm and Jean-Yves
{\sc Thibon}\thanks{Institut Gaspard Monge, Universit\'e de
Marne-la-Vall\'ee, 2 rue de la Butte-Verte, 93166 Noisy-le-Grand cedex,
France}}
\date{}

\begin{document}

\newdimen\Squaresize \Squaresize=14pt
\newdimen\Thickness \Thickness=0.5pt

\def\Square#1{\hbox{\vrule width \Thickness
   \vbox to \Squaresize{\hrule height \Thickness\vss
      \hbox to \Squaresize{\hss#1\hss}
   \vss\hrule height\Thickness}
\unskip\vrule width \Thickness}
\kern-\Thickness}

\def\Vsquare#1{\vbox{\Square{$#1$}}\kern-\Thickness}
\def\blk{\omit\hskip\Squaresize}

\def\young#1{
\vbox{\smallskip\offinterlineskip
\halign{&\Vsquare{##}\cr #1}}}

\maketitle

\begin{abstract}
We introduce a new family of symmetric functions,
which are $q$-analogues of products of
Schur functions, defined
in terms of ribbon tableaux.

These functions can be interpreted in terms of the Fock space
representation $\F$ of $U_q(\slchap_n)$, and are related to
Hall-Littlewood functions via  the geometry
of flag varieties.
We present a series of conjectures, and prove them in  special
cases.
The essential step in proving that these functions are actually
symmetric consists in the calculation of a basis of highest weight
vectors of $\F$ using ribbon tableaux.
\end{abstract}

\footnotesize
\tableofcontents
\normalsize

% ************************************************************************
%                      INTRODUCTION
% ************************************************************************

\section{Introduction}\label{S1}

This article is devoted to the study of a new family of
symmetric functions  $H^{(k)}_\lambda(X;q)$, defined in terms
of certain generalized Young tableaux, called ribbon tableaux,
or rim-hook tableaux \cite{SW}.

These objects, although unfamiliar,
arise naturally in several contexts, and their use is
implicit in many classical algorithms related to the
symmetric groups (see {\it e.g. } Robinson's book \cite{Ro}).
In particular,
they can be applied to the description of the power-sum plethysm
operators $\psi^k : f(\{x_i\})\mapsto f(\{x_i^k\})$ on symmetric functions
\cite{DLT,Mcd}, and this point of view suggests the definition of a natural
$q$-analogue $\psi^k_q$ of $\psi^k$. This $q$-analogue turns out
to make sense when the algebra of symmetric functions is interpreted
as the bosonic Fock space representation of the quantum
affine algebra $U_q(\slchap_k)$. Indeed, one can prove, building
on recent work by Kashiwara, Miwa and Stern \cite{Ste,KMS}, that
the image $\psi^k_q(f)$ of any symmetric function by this operator,
is a highest weight vector for  $U_q(\slchap_k)$.
In particular, the images $\psi^k_q(h_\lambda)$ of products of
complete homogeneous functions have a simple combinatorial
description, and can be used as a convenient basis of highest
weight vectors.

The space of symmetric functions is endowed with a natural
scalar product (the same as in the Fock space interpretation),
and one can consider the adjoint $\varphi^k_q$ of the $q$-plethysm
operator $\psi^k_q$. This operator divides degrees by $k$, and
sends the Schur functions $s_{k\lambda}$ indexed by partitions
of the form $k\lambda=(k\lambda_1,k\lambda_2,\ldots,k\lambda_r)$
onto a new basis, which is essentially the one considered
in this paper. More precisely,
$H^{(k)}_\lambda(X;q^{-2})= \varphi^k_q(s_{k\lambda})$.
It should be said, however, that our original definition was
purely combinatorial,
and that the connection with $U_q(\slchap_k)$
was understood only recently.

The $H$-functions are generalizations of Hall-Littlewood functions.
We prove that for $k$ sufficiently large, $H^{(k)}_\lambda=Q'_\lambda$,
where  $Q'_\lambda$ is the adjoint basis of $P_\lambda$ for the
standard scalar product. Moreover, we conjecture that the
differences $H^{(k+1)}_\lambda-H^{(k)}_\lambda$ are nonnegative
on the Schur basis, \ie that the $H$-functions form a filtration
of the $Q'$-functions. In particular, the coefficients of $H$-functions
on the Schur basis are conjectured to be polynomials with nonnegative
integer coefficients.

The $Q'$-functions are known to be related to a variety of
topics in representation theory \cite{Gr,LLT1,LLT2,LLT4,MS},
algebraic geometry \cite{HSh,HS,La,Lu1,Sh1,Te},
combinatorics \cite{LS1,LS2,Sc}
and mathematical physics \cite{Ki1,KR}. As a general
rule, $q$-analogues related to quantum groups admit interesting
interpretations when the parameter $q$ is specialized to
the cardinality of a finite field, or to a complex root of unity.
The $Q'$-functions are no exception. In the first case,
the coefficients $\tilde K_{\lambda\mu}(q)$ of $\tilde Q'_\mu$
on the Schur basis are character values of the group $GL_n(\GF_q)$ \cite{Lu1},
while in the second one, a factorization property reminescent
of Steinberg's tensor product theorem leads to  combinatorial formulas
for the Schur expansion of certain plethysms,
in particular for $\psi^k(h_\mu)$ \cite{LLT1,LLT2}.

On the basis of extensive numerical computations, we conjecture that
the $H$-functions display the same behaviour with respect to
specializations at roots of unity, giving this time plethysms
$\psi^k(s_\lambda)$ of Schur functions by power sums.
In fact, the $H$-functions were originally defined as $q$-analogues
of products of Schur functions, being the natural generalization
of those introduced in \cite{CL}. A combinatorial description
of general $H$-functions on the Schur basis, similar to the one
given in \cite{CL} in terms of Yamanouchi domino tableaux, would
lead to a refined Littlewood-Richardson rule, compatible with
cyclic symmetrization in the same way as the rule
of \cite{CL} is compatible with symmetrized and antisymmetrized squares.
This means that if one splits a tensor power $V_\lambda^{\otimes k}$
of an irreducible representation $V_\lambda$ of ${\rm U}(n)$ into
eigenspaces $E^{(i)}$ of the cyclic shift operator
$v_1\otimes v_2\otimes\cdots\otimes v_k \mapsto v_2\otimes
v_3\otimes\cdots\otimes
v_k\otimes v_1$,
each $E^{(i)}$ is a representation of ${\rm U}(n)$ whose
spectrum is given, according to the conjectures, by the
coefficient of $q^i$ in the reduction modulo $1-q^k$ of
$H^{(k)}_{\lambda^k}(q)$.

All the conjectures are proved for $k=2$ (domino tableaux)  and
for $k$ sufficiently large (the stable case). The case of domino
tableaux follows from the combinatorial constructions of
\cite{CL} and \cite{KLLT}, while the stable case rely on the
interpretation of Kostka-Foulkes polynomials in terms of characters
of finite linear groups, and as Poincar\'e polynomials of certain
algebraic varieties, in particular on the cell decompositions
of these varieties found by N. Shimomura \cite{Sh1}.

This article is structured as follows.
In Section \ref{HLUV}, we recall some properties of Hall-Littlewood
functions, in particular their interpretation in terms
of affine Hecke algebras and their connection with finite
linear groups.
In Section \ref{S3}, we explain the connexion between plethysm
and Hall-Littlewood functions at roots of unity, and in Section \ref{S4},
we show how to translate these results in terms of ribbon tableaux.
The application of ribbon tableaux to the construction of highest weight
vectors in the Fock representation of $U_q(\slchap_n)$ is presented
in Section \ref{S5}, and connected to a recent construction
of Kashiwara, Miwa and Stern \cite{KMS}. In Section \ref{S6},
we define the $H$-functions, and summarize their known or conjectural
properties. Section \ref{S7} establishes these conjectures for
$H$-functions of level 2, corresponding to domino tableaux.
In Section \ref{S8}, we recall Shimomura's cell decompositions
of unipotent varieties, and show the equivalence of his description
of the Poincar\'e polynomials with a variant  needed in the sequel.
{}From this, we deduce
that the $H$-functions of sufficiently
large level are equal to Hall-Littlewood functions, which is also
sufficient to prove all the conjectured properties in this case.

\section{Hall-Littlewood functions}\label{HLUV}

Our notations for symmetric functions will be essentially
those of the book \cite{Mcd},
to which the reader is referred for more details.

The original definition of Hall-Littlewood functions can be
reformulated in terms of an action
of the affine Hecke algebra $\widehat{H}_N(q)$ of type $A_{N-1}$ on
the ring $\C[x_1^{\pm 1},\ldots, x_N^{\pm 1}]$ \cite{DKLLST}.

The affine Hecke algebra $\widehat{H}_N(q)$ is generated by
$T_i$, $i=1,\ldots,N-1$ and $y_i^{\pm 1}$, $i=1,\ldots,N$,
with relations
\begin{equation}
\left\{\matrix{
T_i^2  =  (q-1)T_i +q \cr
T_i T_{i+1} T_i  =  T_{i+1}T_i T_{i+1} \cr
T_i T_j = T_j T_i \ \ (|j-i|>1) \cr }\right.
\qquad\quad
\left\{\matrix{
y_i y_j=y_jy_i \cr
y_jT_i=T_iy_j \quad j\not= i,i+1\cr
y_jT_j=T_j y_{j+1} -(q-1)y_{j+1} \cr
y_{j+1}T_j = T_j y_j +(q-1) y_{j+1} \cr}\right.
\end{equation}

If $\sigma=\sigma_{i_1}\sigma_{i_2}\cdots\sigma_{i_r}$ is a reduced
decomposition of a permutation $\sigma\in\S_N$, where $\sigma_i=(i,i+1)$,
one sets as usual
$T_\sigma=T_{i_1}T_{i_2}\ldots T_{i_r}$, the result being
independent of the reduced decomposition.

Let $\displaystyle \Delta_N(q)=\prod_{1\le i<j \le N}(qx_i-x_j)$.
Then, on one hand, the Hall-Littlewood polynomial
$Q_\lambda (x_1,\ldots,x_N ; q)$ indexed by a partition $\lambda$
of length $\le N$ is defined by \cite{Li1}
\begin{equation}
Q_\lambda = {(1-q)^{\ell(\lambda)}\over [m_0]_q !}
\sum_{\sigma\in\S_N} \sigma
\left( x^\lambda {\Delta_N(q)\over \Delta_N(1)} \right)
\end{equation}
where $m_0=N-\ell(\lambda)$ and the $q$-integers are here defined
by
$[n]_q=(1-q^n)/(1-q)$.

On the other hand, $\widehat{H}_N(q)$ acts on  $\C[x_1^{\pm 1},\ldots, x_N^{\pm
1}]$
by $y_i(f)=x_if$ and $T_i=(q-1)\pi_i+\sigma_i$, where
$\pi_i$ is the isobaric divided difference operator
$$
\pi_i(f) = {x_if -x_{i+1}\sigma_i(f)\over x_i-x_{i+1}} \ ,
$$
and it is shown in \cite{DKLLST} that if one defines the $q$-symmetrizing
operator $S^{(N)}\in \widehat{H}_N(q)$ by
\begin{equation}
S^{(N)} = \sum_{\sigma\in\S_N} T_\sigma \ ,
\end{equation}
then
\begin{equation}
Q_\lambda(x_1,\ldots,x_N ; q) = {(1-q)^{\ell(\lambda)}\over [m_0]_q!} S^{(N)}
(x^\lambda) \ .
\end{equation}

The normalization factor $1/[m_0]_q!$ is here to ensure stability
with respect to the adjunction of variables, and if we denote
by $X$ the infinite set $X=\{x_1,x_2,\ldots,\, \}$ then
$Q_\lambda(X;q) = \lim_{N\rightarrow\infty} Q_\lambda(x_1,\ldots,x_N;q)$.

The $P$-functions are defined by
$$
P_\lambda(X;q) =
{1\over (1-q)^{\ell(\lambda)}[m_1]_q!\cdots [m_n]_q!}Q_\lambda(X;q)
$$
where $m_i$ is the multiplicity of the part $i$ in $\lambda$.

We consider these functions as elements of the algebra $\sym=\sym(X)$
of symmetric functions with coefficients in $\C(q)$.
In this paper, the scalar product $\<\, ,\,\>$ on $\sym$
will always be the standard one, for which the Schur functions $s_\lambda$
form an orthonormal basis.

We denote
by $(Q'_\mu(X;q))$ the adjoint basis of $P_\lambda(X;q)$ for this
scalar product. It is easy to see that $Q'_\mu(X;q)$ is
the image of  $Q_\mu(X;q)$
by the ring homomorphism $p_k\mapsto (1-q^k)^{-1}p_k$
(in $\lambda$-ring notation, $Q'_\mu(X;q)=Q(X/(1-q);q)$).
In the Schur basis,
\begin{equation}
Q'_\mu(X;q)=\sum_\lambda K_{\lambda\mu}(q)s_\lambda(X)
\end{equation}
where the $K_{\lambda\mu}(q)$ are the Kostka-Foulkes polynomials.
The polynomial $K_{\lambda\mu}(q)$ is the generating function
of a statistic $c$ called ${\it charge}$ on the set $\tab(\lambda,\mu)$ of
Young
tableaux of shape $\lambda$ and weight $\mu$ (\cite{LS2,Sc}, see also
\cite{La,Mcd})
\begin{equation}
K_{\lambda\mu}(q)=\sum_{\t\in\tab(\lambda,\mu)}q^{c(\t)} \ .
\end{equation}
We shall also need the $\tQ$-functions, defined by
\begin{equation}
\tQ_\mu(X;q) = \sum_\lambda \tK_{\lambda\mu}(q)s_\lambda(X)
=q^{n(\mu)}Q'_\mu(X;q^{-1}) \ .
\end{equation}
The polynomial $\tK_{\lambda\mu}(q)$ is the generating function
of the complementary statistic $\cc(\t) = n(\mu)-c(\t)$, which is called
{\it cocharge}. The operation of {\it cyclage} endows $\tab(\lambda,\mu)$
with the structure of a rank poset, in which the rank of a tableau is
equal to its cocharge (see \cite{La}).

When the parameter $q$ is interpreted as the cardinality of a
finite field $\GF_q$, it is known that $\tK_{\lambda\mu}(q)$ is equal
to the value $\chi^\lambda(u)$ of the unipotent character $\chi^\lambda$
of $G=GL_n(\GF_q)$ on a unipotent element $u$ with Jordan canonical form
specified by the partition $\mu$ (see \cite{Lu2}).

In this specialization, the coefficients
\begin{equation}
\tG_{\nu\mu}(q) = \<h_\nu\, ,\,\tQ_\mu\>
\end{equation}
of the $\tQ$-functions on the basis of monomial symmetric functions are
also the values of certain characters of $G$ on unipotent classes. Let
$\P_\nu$ denote a parabolic subgroup of type $\nu$ of $G$, for example
the group of upper block triangular matrices with diagonal blocks of sizes
$\nu_1,\ldots,\nu_r$, and consider the permutation representation of
$G$ over $\C[G/\P_\nu]$. The value $\xi^\nu(g)$ of the character $\xi^\nu$ of
this
representation on an element $g\in G$ is equal to the number of fixed
points of $g$ on $ G/\P_\nu$. Then, it can be shown that, for a unipotent
$u$ of type $\mu$,
\begin{equation}
\xi^\nu(u)=\tG_{\nu\mu}(q) \ .
\end{equation}
The factor set $G/\P_\nu$ can be identified with the variety $\F_\nu$ of
$\nu$-flags in $V=\GF_q^n$
$$
V_{\nu_1}\subset V_{\nu_1+\nu_2}\subset\ldots\subset V_{\nu_1+\ldots \nu_r}=V
$$
where $\dim V_i = i$. Thus, $\tG_{\nu\mu}(q)$ is equal to the number
of $\GF_q$-rational points of the algebraic variety $\F_\nu^u$ of fixed
points of $u$ in $\F_\nu$.

\section{Specializations at roots of unity}\label{S3}

As recalled in the preceding section, the Hall-Littlewood functions with
parameter specialized to the cardinality $q$ of a finite field $\GF_q$
provide information about the complex
characters of the linear group $GL(n,\GF_q)$
over this field. It turns out that when the parameter is specialized
to a complex root of unity, one obtains information about representations of
$GL(n,\C)$ (or ${\rm U}(n)$,
that is, a combinatorial decomposition of certain plethysms
\cite{LLT1,LLT2}. We give now a brief review of these results.

The first one is a factorization property of the functions $Q'_\lambda(X,q)$
when $q$ is specialized to a primitive root of unity. This is to be seen as
a generalization of the fact that when $q$ is specialized to $1$ the function
$Q'_{\lambda}(X;q)$ reduces to $h_{\lambda}(X) = \prod_i h_{\lambda_i}(X)$.

\begin{theorem}{\rm \cite{LLT1} }\label{THLLT1}
Let $\lambda = (1^{m_1}2^{m_2} \ldots n^{m_n})$ be a
partition written multiplicatively.
Set $m_i=kq_i+r_i$ with $0\le r_i<k$, and
$\mu=(1^{r_1}2^{r_2} \ldots n^{r_n})$. Then, $\zeta$ being a primitive
$k$-th root of unity,
\begin{equation}\label{FACT}
Q'_{\lambda}(X;\zeta)=Q'_\mu(X;\zeta)
     \prod_{i\ge 1}\bigl[ Q'_{(i^k)}(X;\zeta)\bigr]^{q_i} \ .
\end{equation}
\end{theorem}

The functions $Q'_{(i^k)}(X;\zeta)$ appearing in the right-hand side of
(\ref{FACT})
can be expressed as plethysms.

\begin{theorem}{\rm \cite{LLT1} }\label{THLLT2}
Let $p_k\circ h_n$ denote the plethysm of the
complete function $h_n$ by the power-sum $p_k$, which is defined
by the generating series
$$ \sum_n  p_k\circ h_n(X)\,z^n=\prod_{x\in X}(1-zx^k)^{-1}\ . $$
Then, if $\zeta$ is as above a primitive $k$-th root of unity, one has
$$Q'_{(n^k)}(X;\zeta)=(-1)^{(k-1)n}p_k\circ h_n(X). $$
\end{theorem}

For example, with $k=3$ $(\zeta = e^{2i\pi /3})$, we have
$$
Q'_{444433311}(X;\zeta)=Q'_{411}(X;\zeta)\,Q'_{333}(X;\zeta)\,Q'_{444}(X;\zeta)
=Q'_{411}(X;\zeta)\,p_4\circ h_{43} \ .
$$

Let $V$ be a polynomial representation of $GL(n,\C)$, with
character the symmetric function $f$. Let $\gamma$ be
the cyclic shift operator on $V^{\otimes k}$, that is,
$$
\gamma (v_1\otimes v_2\otimes \cdots\otimes v_k)
=
v_2\otimes v_3\otimes\cdots\otimes v_k\otimes v_1 \ .
$$
Let $\zeta =\exp(2{\rm i}\pi/k)$, and denote by $E^{(r)}$
the eigenspace of $\gamma$ in $V^{\otimes k}$ associated
with the eigenvalue $\zeta^r$.
As $\gamma$ commutes with the action of $GL(n)$, these eigenspaces
are representations of $GL(n)$, and their characters are given
by the plethysms
$\ell_k^{(r)}\circ f$
of the character $f$ of $V$ by certain symmetric functions
$\ell_k^{(r)}$ that we shall now describe.

For $k,n\in {\bf N}$, the \it Ramanujan \rm or \it Von Sterneck
\it sum \rm $c(k,n)$ (also denoted $\Phi(k,n)$) is the sum of the
$k$-th powers of the \it primitive \rm $n$-th roots of unity.  Its value
is given by \it H\"older's formula\rm: if $(k,n)=d$ and $n=md$, then
$c(k,n)=\mu(m)\phi(n)/\phi(m)$, where $\mu$ is the Moebius function and
$\phi$ is the Euler totient function (see {\it e.g.}  \cite{NV}).

The  symmetric functions $\ell_k^{(r)}$ are given by the formula
\begin{equation}
\ell^{(r)}_k={{1}\over{k}}\sum_{d\mid k}c(r,d)p_d^{k/d} \ .
\end{equation}
These functions are the Frobenius characteristics
of the representations of the symmetric group induced by irreducible
representations of a transitive  cyclic subgroup \cite{Fo}.
A combinatorial interpretation of the multiplicity $\<s_\lambda\, ,\,
\ell^{(k)}_n\>$
has been given by Kraskiewicz and Weyman \cite{KW}. This result is equivalent
to the congruence
$$
Q'_{1^n}(X;q) \equiv \sum_{0\le k \le n-1} q^k \ell_n^{(k)} \ (\mod 1-q^n) \ .
$$
Another  proof can be found in \cite{De}.
Now, if $V$ is a product of exterior powers of the fundamental
representation $\C^n$
$$
V=\Lambda^{\nu_1}\C^n\otimes \Lambda^{\nu_2}\C^n\otimes\cdots\otimes
\Lambda^{\nu_m}\C^n
$$
the character of $E^{(r)}$ is $\ell_k^{(r)}\circ e_\nu$, and
similarly if $V$ is a product of symmetric powers with
character $h_\nu$, the character of $E^{(r)}$ is
$\ell_k^{(r)}\circ h_\nu$.

Given two partitions $\lambda$ and $\mu$, we denote by $\lambda \vee \mu$
the partition obtained by reordering the concatenation of $\lambda$
and $\nu$, {\it e.g.}  $(2,\,2,\,1)\vee (5,\,2,\,1,\,1) = (5,\,2^3,\,1^3)$.
We write $\mu^k=\mu\vee \mu\vee \cdots\vee \mu$
($k$ factors). If $\mu = (\mu_1,\,\ldots ,\,\mu_r)$, we set
$k\mu = (k\mu_1,\,\ldots ,\,k\mu_r)$.

Taking into account Theorems \ref{THLLT1} and \ref{THLLT2}
and following the method of \cite{De}, one arrives at the
following combinatorial formula for the decomposition
of $E^{(r)}$ into irreducibles:

\begin{theorem}{\rm \cite{LLT2} }
 Let $e_i$ be the $\ i$-th elementary symmetric
function, and for $\lambda=(\lambda_1,\ldots,\lambda_m)$, $e_\lambda=
e_{\lambda_1}\cdots e_{\lambda_r}$.
Then, the multiplicity $\<s_\mu\, ,\,\ell^{(r)}_k\circ e_\lambda\>$
of the Schur function $s_\mu$ in the plethysm $\ell^{(r)}_k\circ e_\lambda$
is equal to the number of Young tableaux of shape $\mu'$ (conjugate
partition) and weight $\lambda^k$ whose charge is congruent to $r$ modulo $k$.

This gives as well the plethysms with product of complete functions, since
$$
\<s_{\mu'}\, ,\, \ell_k^{(r)}\circ e_\lambda \> =
\left\{ \matrix{
\<s_\mu\, ,\, \ell_k^{(r)}\circ h_\lambda\> & \mbox{if $|\lambda|$ is even}\cr
\<s_\mu\, ,\, \tilde\ell_k^{(r)}\circ h_\lambda\> & \mbox{if $|\lambda|$ is
odd}\cr
}\right.
$$
where $\tilde\ell_k^{(r)}=\omega(\ell_k^{(r)})=\ell_k^{(s)}$ with
$s=k(k-1)/2-r$.
\end{theorem}

For example, with $k=4$, $r=2$ and $\lambda=(2)$,
$$\ell^{(2)}_4\circ e_2 = s_{431} +s_{422} +s_{41111}+2s_{3311}+2s_{3221}$$
$$+2s_{32111} +s_{2222}+s_{22211}+2s_{221111}+s_{2111111} \ .$$
The coefficient $\<s_{32111}\, ,\, \ell^{(2)}_4\circ e_2\>=2$
is equal to the number of tableaux of shape $(3,\,2,\,1,\,1,\,1)'=(5,\,2,\,1)$,
weight $(2,\,2,\,2,\,2)$
and charge $\equiv 2\ (\mod 4)$. The two tableaux satisfying
these conditions are

$$
\young{3\cr 2&4\cr 1&1&2&3&4\cr}\qquad\qquad
\young{4\cr 2&3 \cr 1&1&2&3&4\cr}
$$

\bigskip\noindent
which both have charge equal to $6$. \par
Similarly,  $\<s_{732}\, ,\, \ell^{(2)}_4\circ e_{21}\>=5$
is the number of tableaux with shape $(3,\,3,\,2,\,1,\,1,\,1,\,1)$, weight
$(2,\,2,\,2,\,2,\,1,\,1,\,1,\,1)$ and
charge $\equiv 2\ (\mod 4)$.

Another  combinatorial formulation of Theorems \ref{THLLT1} and
\ref{THLLT2} can be presented by means of the notion of
{\it ribbon tableau}, which will also provide the clue for their
generalization.

\section{Ribbon tableaux}\label{S4}

To a partition $\lambda$ is associated a $k$-core $\lambda_{(k)}$
and a $k$-quotient $\lambda^{(k)}$
\cite{JK}. The $k$-core is the unique partition obtained by successively
removing
$k$-ribbons (or skew hooks) from $\lambda$. The different possible ways of
doing so can be distinguished from one another by labelling $1$ the
last ribbon removed, $2$ the penultimate, and so on. Thus Figure~\ref{TRIBB}
shows two different ways of reaching the $3$-core $\lambda_{(3)}=(2,\,1^2)$
of $\lambda = (8,\, 7^2,\,4,\,1^5)$. These pictures represent two $3$-ribbon
tableaux $T_1,\,T_2$ of shape $\lambda/\lambda_{(3)}$ and weight $\mu = (1^9)$.

\begin{figure}[h]\label{FIG1}
\setlength{\unitlength}{0.25pt}
\centerline{
\begin{picture}(1400,500)(0,0)
\put(0,10){$T_1 =$}
\put(120,0){\begin{picture}(400,450)
\put(0,0){\framebox(50,50){}}
\put(50,0){\framebox(50,50){}}
\put(0,50){\framebox(50,50){}}
\put(0,100){\framebox(50,50){}}
\put(100,0){\framebox(150,50){}}
\put(0,150){\framebox(50,150){}}
\put(0,300){\framebox(50,150){}}
\put(50,150){\framebox(150,50){}}
\put(200,100){\framebox(150,50){}}
\put(150,50){\line(0,1){50}}
\put(150,100){\line(-1,0){50}}
\put(100,100){\line(0,1){50}}
\put(200,50){\line(0,1){50}}
\put(300,0){\line(0,1){100}}
\put(350,50){\line(0,1){50}}
\put(350,50){\line(1,0){50}}
\put(400,0){\line(0,1){50}}
\put(250,0){\line(1,0){150}}
\put(165,10){1}
\put(65,60){2}
\put(115,110){3}
\put(15,160){4}
\put(215,60){5}
\put(165,160){6}
\put(315,10){7}
\put(15,310){8}
\put(265,110){9}
\end{picture}}
\put(880,10){$T_2 = $}
\put(1000,0){\begin{picture}(400,450)
\put(0,0){\framebox(50,50){}}
\put(50,0){\framebox(50,50){}}
\put(0,50){\framebox(50,50){}}
\put(0,100){\framebox(50,50){}}
\put(250,0){\framebox(150,50){}}
\put(0,150){\framebox(50,150){}}
\put(0,300){\framebox(50,150){}}
\put(50,150){\framebox(150,50){}}
\put(200,100){\framebox(150,50){}}
\put(50,100){\framebox(150,50){}}
\put(200,50){\framebox(150,50){}}
\put(150,0){\line(0,1){100}}
\put(250,0){\line(0,1){50}}
\put(100,0){\line(1,0){150}}
\put(165,10){4}
\put(65,60){2}
\put(115,110){6}
\put(15,160){1}
\put(215,60){7}
\put(165,160){8}
\put(315,10){5}
\put(15,310){3}
\put(265,110){9}
\end{picture}}
\end{picture}}
\caption{\label{TRIBB}}
\end{figure}

To define  $k$-ribbon tableaux of general weight and shape, we need some
terminology.
The {\it initial cell} of a $k$-ribbon $R$ is its rightmost and bottommost
cell. Let $\theta = \beta/\alpha$ be a skew shape, and set
$\alpha_+ = (\beta_1)\vee \alpha$, so that $\alpha_+/\alpha$ is the horizontal
strip made of the bottom cells of the columns of $\theta$. We say that $\theta$
is a {\it horizontal $k$-ribbon strip} of weight $m$, if it can be tiled by
$m$ $k$-ribbons the initial cells of which lie in $\alpha_+/\alpha$. (One can
check that if such a tiling exists, it is unique).

Now, a {\it $k$-ribbon tableau} $T$ of shape $\lambda/\nu$ and weight
$\mu=(\mu_1,\,\ldots ,\,\mu_r)$ is defined as a chain of partitions
$$
\nu=\alpha^0\subset \alpha^1 \subset \cdots \subset \alpha^r=\lambda
$$
such that $\alpha^i/\alpha^{i-1}$ is a horizontal $k$-ribbon strip of weight
$\mu_i$. Graphically, $T$ may be described by numbering each $k$-ribbon of
$\alpha^i/\alpha^{i-1}$ with the number $i$. We denote by
$\tab_k(\lambda/\nu,\,\mu)$ the
set of $k$-ribbon tableaux of shape $\lambda/\nu$ and weight $\mu$, and we set
$$
K_{\lambda/\nu,\,\mu}^{(k)} = |\tab_k(\lambda/\nu,\,\mu)| \ .
$$
Finally we recall the definition of the $k$-sign $\epsilon_k(\lambda/\nu)$.
Define
the sign of a ribbon $R$ as $(-1)^{h-1}$, where $h$ is the height of $R$. The
$k$-sign $\epsilon_k(\lambda/\nu)$ is the product of the signs of all the
ribbons
of a $k$-ribbon tableau of shape $\lambda/\nu$ (this does not depend on the
particular tableau chosen, but only on the shape).

The origin of these combinatorial definitions is best understood by analyzing
carefully the operation of multiplying a Schur function $s_\nu$ by a plethysm
of the form $\psi^k(h_\mu)=p_k \circ h_\mu$.
Equivalently, thanks to the involution $\omega$,
one may rather consider a product of the type $s_\nu \, [p_k\circ e_\mu]$. To
this end,
since
$$
p_k\circ e_\mu = (e_{\mu_1}\circ p_k)\, \cdots (e_{\mu_n}\circ p_k)
= m_{k^{\mu_1}}\cdots m_{k^{\mu_n}}
$$
one needs only to apply repeatedly the following multiplication rule due
to Muir \cite{Mu} (see also \cite{Li3}):
$$
s_\nu \, m_\alpha = \sum_\beta s_{\nu + \beta} \ ,
$$
sum over all distinct permutations $\beta$ of
$(\alpha_1,\,\alpha_2,\,\ldots ,\, \alpha_n,\,0,\,\ldots \, )$.
Here  the Schur functions
$s_{\nu + \beta}$ are not necessarily indexed by partitions and have therefore
to be put in standard form, this reduction yielding only a finite number of
nonzero
summands. For example,
$$
s_{31}\,m_3 = s_{61} + s_{313} + s_{31003} = s_{61} - s_{322} + s_{314} \ .
$$
Other terms such as $s_{34}$ or $s_{3103}$
reduce to $0$. It is easy to deduce
from this rule that the multiplicity
$$
\< s_\nu \, m_{k^{\mu_i}} \, , \, s_\lambda \>
$$
is nonzero iff $\lambda '/\mu '$ is a horizontal $k$-ribbon strip of weight
$\mu_i$, in which case it is equal to $\epsilon_k(\lambda/\nu)$. Hence,
applying
$\omega$ we arrive at the expansion
\begin{equation}\label{plethrub}
s_\nu \, [p_k\circ h_\mu] = \sum_\lambda \epsilon_k(\lambda/\mu) \,
K_{\lambda/\nu \,, \mu}^{(k)} \, s_\lambda
\end{equation}
from which we deduce by \ref{THLLT1}, \ref{THLLT2} that
$$
K_{\lambda \, \mu}^{(k)} = (-1)^{(k-1)|\mu|} \, \epsilon_k(\lambda)
\, K_{\lambda \, \mu^k}(\zeta)
$$
and more generally, defining as in \cite{KR} the skew Kostka-Foulkes polynomial
$K_{\lambda/\nu \,, \alpha}(q)$ by
$$
K_{\lambda/\nu \,, \alpha}(q) = \< s_{\lambda/\nu} \, , \, Q'_\alpha(q) \>
$$
(or as the generating functions of the charge statistic on skew
tableaux of shape $\lambda/\mu$ and weight $\alpha$),
we can write
$$
K_{\lambda/\nu \,, \mu}^{(k)} = (-1)^{(k-1)|\mu|} \, \epsilon_k(\lambda/\nu)
\, K_{\lambda/\nu \,, \mu^k}(\zeta) \ .
$$

It turns out that enumerating $k$-ribbon tableaux is equivalent to enumerating
$k$-uples of ordinary Young tableaux, as shown by the correspondence to be
described
now. This bijection was first introduced by Stanton and White \cite{SW} in the
case
of ribbon tableaux of right shape $\lambda$ (without $k$-core) and standard
weight $\mu = (1^n)$ (see also \cite{FS}). We need some additional definitions.

Let $R$ be a $k$-ribbon of a $k$-ribbon tableau. $R$ contains a unique cell
with coordinates $(x,\,y)$
such that $y-x\equiv 0 \ (\mod k)$. We decide to write in this cell the number
attached to
$R$, and we define the {\it type} $i\in \{0,\,1,\,\ldots ,\,k-1\}$
of $R$ as the distance between this cell and the initial cell of $R$.
For example, the $3$-ribbons of
$T_1$ are divided up into three classes (see Fig. \ref{TRIBB}):
\begin{itemize}
\item 4, 6, 8, of type 0;
\item 1, 2, 7, 9, of type 1;
\item 3, 5, of type 2.
\end{itemize}
Define the {\it diagonals} of a $k$-ribbon tableau as the sequences of integers
read along the straight lines $D_i \, : \, y-x = ki $.
Thus $T_1$ has the sequence of diagonals
$$((8),\,(4),\, (2,\,3,\,6),\,
(1,\,5,\,9),\,(7))\ .$$
This definition applies in particular to $1$-ribbon tableaux, \ie ordinary
Young tableaux. It is obvious that a Young tableau is uniquely determined
by its sequence of diagonals. Hence, we can associate to a given $k$-ribbon
tableau $T$ of
shape $\lambda/\nu$ a $k$-uple $(t_0,\,t_1,\,\ldots ,\,t_{k-1})$ of
Young tableaux defined as follows; the diagonals of $t_i$ are obtained by
erasing
in the diagonals of $T$ the labels of all the ribbons of type $\not = i$.
For instance, if $T=T_1$ the first ribbon tableau of Figure~\ref{TRIBB}, the
sequence of
diagonals of $t_1$ is $\left((2),\,(1,\,9),\,(7)\right)$, and

\setlength{\unitlength}{0.25pt}
\centerline{
\begin{picture}(300,150)
\put(100,0){\begin{picture}(100,100)
\put(0,0){\framebox(50,50){1}}
\put(50,0){\framebox(50,50){7}}
\put(0,50){\framebox(50,50){2}}
\put(50,50){\framebox(50,50){9}}
\end{picture}}
\put(0,20){$t_1 =$}
\end{picture}
}

\bigskip
\noindent
The complete triple $(t_0,\,t_1,\,t_2)$ of Young tableaux associated to $T_1$
is

\setlength{\unitlength}{0.25pt}
\centerline{
\begin{picture}(700,150)
\put(150,0){\begin{picture}(100,100)
\put(0,0){\framebox(50,50){4}}
\put(50,0){\framebox(50,50){6}}
\put(0,50){\framebox(50,50){8}}
\end{picture}}
\put(320,0){\begin{picture}(100,100)
\put(0,0){\framebox(50,50){1}}
\put(50,0){\framebox(50,50){7}}
\put(0,50){\framebox(50,50){2}}
\put(50,50){\framebox(50,50){9}}
\end{picture}}
\put(490,0){\begin{picture}(100,100)
\put(0,0){\framebox(50,50){3}}
\put(50,0){\framebox(50,50){5}}
\end{picture}}
\put(0,40){$\tau^1 =$}
\put(110,40){$\Big( $}
\put(270,40){,}
\put(440,40){,}
\put(610,40){$\Big) $}
\end{picture}}

\bigskip
\noindent
whereas that corresponding to $T_2$ is

\setlength{\unitlength}{0.25pt}
\centerline{
\begin{picture}(700,150)
\put(150,0){\begin{picture}(100,100)
\put(0,0){\framebox(50,50){1}}
\put(50,0){\framebox(50,50){8}}
\put(0,50){\framebox(50,50){3}}
\end{picture}}
\put(320,0){\begin{picture}(100,100)
\put(0,0){\framebox(50,50){4}}
\put(50,0){\framebox(50,50){5}}
\put(0,50){\framebox(50,50){6}}
\put(50,50){\framebox(50,50){9}}
\end{picture}}
\put(490,0){\begin{picture}(100,100)
\put(0,0){\framebox(50,50){2}}
\put(50,0){\framebox(50,50){7}}
\end{picture}}
\put(0,40){$\tau^2 =$}
\put(110,40){$\Big( $}
\put(270,40){,}
\put(440,40){,}
\put(610,40){$\Big) $}
\end{picture}}

\bigskip
\noindent
One can show that if $\nu =\lambda_{(k)}$, the $k$-core of $\lambda$,
the $k$-uple of shapes $(\lambda^0,\, \lambda^1,\, \ldots,\,\lambda^{k-1})$
of $(t_0,\,t_1,\,\ldots ,\,t_{k-1})$ depends only on the shape $\lambda$ of
$T$, and is equal to the $k$-quotient $\lambda^{(k)}$ of $\lambda$.
Moreover the correspondence $T \longrightarrow (t_0,\,t_1,\,\ldots ,\,t_{k-1})$
establishes a
bijection between the set of $k$-ribbon tableaux of shape
$\lambda/\lambda_{(k)}$
and weight $\mu$, and the set of $k$-uples of Young tableaux of
shapes $(\lambda^0,\,\ldots ,\,\lambda^{k-1})$ and weights $(\mu^0,\, \ldots
,\,\mu^{k-1})$
with $\mu_i = \sum_j \mu^j_i$.
(See \cite{SW} or \cite{FS}
for a proof in the case when $\lambda_{(k)}= (0)$ and $\mu = (1^n)$).

For example, keeping $\lambda = (8,\, 7^2,\,4,\,1^5)$, the triple

\centerline{
\begin{picture}(700,150)                \setlength{\unitlength}{0.25pt}
\put(150,0){\begin{picture}(100,100)
\put(0,0){\framebox(50,50){3}}
\put(50,0){\framebox(50,50){3}}
\put(0,50){\framebox(50,50){4}}
\end{picture}}
\put(320,0){\begin{picture}(100,100)
\put(0,0){\framebox(50,50){1}}
\put(50,0){\framebox(50,50){3}}
\put(0,50){\framebox(50,50){2}}
\put(50,50){\framebox(50,50){4}}
\end{picture}}
\put(490,0){\begin{picture}(100,100)
\put(0,0){\framebox(50,50){2}}
\put(50,0){\framebox(50,50){3}}
\end{picture}}
\put(0,40){$\tau =$}
\put(110,40){$\Big( $}
\put(270,40){,}
\put(440,40){,}
\put(610,40){$\Big) $}
\end{picture}}

\bigskip
\noindent
with weights $\left(
(0,\,0,\,2,\,1),\,(1,\,1,\,1,\,1),\,(0,\,1,\,1,\,0)\right)$
corresponds to the 3-ribbon tableau

\centerline{
\begin{picture}(550,500)(0,0)              \setlength{\unitlength}{0.25pt}
\put(0,10){$T =$}
\put(120,0){\begin{picture}(400,450)
\put(0,0){\framebox(50,50){}}
\put(50,0){\framebox(50,50){}}
\put(0,50){\framebox(50,50){}}
\put(0,100){\framebox(50,50){}}
\put(100,0){\framebox(150,50){}}
\put(0,150){\framebox(50,150){}}
\put(0,300){\framebox(50,150){}}
\put(50,150){\framebox(150,50){}}
\put(200,100){\framebox(150,50){}}
\put(150,50){\line(0,1){50}}
\put(150,100){\line(-1,0){50}}
\put(100,100){\line(0,1){50}}
\put(200,50){\line(0,1){50}}
\put(300,0){\line(0,1){100}}
\put(350,50){\line(0,1){50}}
\put(350,50){\line(1,0){50}}
\put(400,0){\line(0,1){50}}
\put(250,0){\line(1,0){150}}
\put(165,10){1}
\put(65,60){2}
\put(115,110){2}
\put(15,160){3}
\put(215,60){3}
\put(165,160){3}
\put(315,10){3}
\put(15,310){4}
\put(265,110){4}
\end{picture}}
\end{picture}}

\bigskip
\noindent
of weight $\mu=(1,\,2,\,4,\,2)$.

As before, the significance of this combinatorial construction becomes clearer
once interpreted in terms of symmetric functions. Recall the definition of
$\phi_k$, the adjoint of the linear operator
$\psi^k:\ F\mapsto p_k\circ F$ acting on the space of symmetric functions.
In other words, $\phi_k$ is characterized by
$$
\< \phi_k(F) \, , \, G \> = \< F \, , \, p_k \circ G \> \ , \ \ \  F,\,G \in
\sym \ .
$$
Littlewood has shown \cite{Li3} that if $\lambda$ is a partition whose
$k$-core $\lambda_{(k)}$ is null, then
\begin{equation}\label{PHIQUOT}
\phi_k(s_\lambda) = \epsilon_k(\lambda) \, s_{\lambda^0} \, s_{\lambda^1} \,
\cdots
\, s_{\lambda^{k-1}}
\end{equation}
where $\lambda^{(k)} = (\lambda^0 ,\, \ldots \, ,\, \lambda^{k-1} )$ is the
$k$-quotient. Therefore,
$$
K_{\lambda \, \mu}^{(k)} = \epsilon_k(\lambda) \, \< p_k \circ h_\mu \, , \,
s_\lambda \>
= \epsilon_k(\lambda) \, \< \phi_k(s_\lambda) \, , \, h_\mu \>
= \< s_{\lambda^0} \, s_{\lambda^1} \, \cdots
\, s_{\lambda^{k-1}} \, , \, h_\mu \>
$$
is the multiplicity of the weight $\mu$ in the product of Schur functions
$s_{\lambda^0} \, \cdots \, s_{\lambda^{k-1}}$, that is, is equal to the number
of
$k$-uples of Young tableaux of shapes $(\lambda^0,\,\ldots ,\,\lambda^{k-1})$
and
weights $(\mu^0,\, \ldots ,\,\mu^{k-1})$ with $\mu_i = \sum_j \mu^j_i$. Thus,
the
bijection described above gives a combinatorial proof of (\ref{PHIQUOT}).

More generally, if $\lambda$ is replaced by a skew partition $\lambda/\nu$,
(\ref{PHIQUOT}) becomes \cite{KSW}
$$
\phi_k(s_{\lambda/\nu}) = \epsilon_k(\lambda/\nu) \, s_{\lambda^0/\nu^0} \,
s_{\lambda^1/\nu^1} \, \cdots \, s_{\lambda^{k-1}/\nu^{k-1}}
$$
if $\lambda_{(k)} = \nu_{(k)}$, and $0$ otherwise. This can also be deduced
from the
previous combinatorial correspondence, but we shall not go into further
details.

Returning to Kostka polynomials, we may summarize this discussion by stating
Theorems~\ref{THLLT1} and \ref{THLLT2} in the following way:

\begin{theorem}
Let $\lambda$ and $\nu$ be partitions and set $\nu = \mu^k \vee \alpha$ with
$m_i(\alpha) < k$. Denoting by $\zeta$ a primitive $k$th root of unity, one
has
\begin{equation}\label{KOSTROOT}
K_{\lambda,\,\nu}(\zeta) = (-1)^{(k-1)|\mu|}
\sum_\beta \epsilon_k(\lambda/\beta)\, K_{\lambda/\beta,\,\mu}^{(k)}\,
K_{\beta,\,\alpha}(\zeta) \ .
\end{equation}
\end{theorem}

\begin{example}{\rm We take $\lambda = (4^2,\,3)$, $\nu = (2^2,\,1^7)$ and
$k=3$ $(\zeta = e^{2i\pi/3})$. In this case,
$\nu= \mu^k \vee \alpha$ with $\mu = (1^2)$ and
$\alpha = (2^2,\,1)$. The summands of (\ref{KOSTROOT}) are parametrized
by the $3$-ribbon tableaux of external shape $\lambda$ and weight $\mu$.
Here we have three such tableaux:

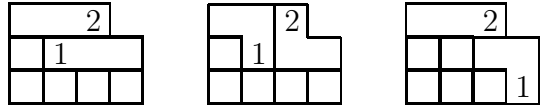
\centerline{
\begin{picture}(800,250)              \setlength{\unitlength}{0.25pt}
\put(0,0){\begin{picture}(200,150)
\put(0,0){\framebox(50,50){}}
\put(50,0){\framebox(50,50){}}
\put(100,0){\framebox(50,50){}}
\put(150,0){\framebox(50,50){}}
\put(0,50){\framebox(50,50){}}
\put(0,0){\line(1,0){200}}
\put(0,100){\line(1,0){200}}
\put(0,0){\line(0,1){100}}
\put(200,0){\line(0,1){100}}
\put(50,50){\framebox(150,50){}}
\put(0,100){\framebox(150,50){}}
\put(65,60){1}
\put(115,110){2}
\end{picture}}
\put(300,0){\begin{picture}(200,150)
\put(0,0){\framebox(50,50){}}
\put(50,0){\framebox(50,50){}}
\put(100,0){\framebox(50,50){}}
\put(150,0){\framebox(50,50){}}
\put(0,50){\framebox(50,50){}}
\put(0,0){\line(1,0){200}}
\put(0,150){\line(1,0){150}}
\put(0,100){\line(1,0){50}}
\put(50,50){\line(0,1){50}}
\put(50,50){\line(1,0){150}}
\put(100,50){\line(0,1){100}}
\put(0,0){\line(0,1){150}}
\put(200,0){\line(0,1){100}}
\put(150,150){\line(0,-1){50}}
\put(150,100){\line(1,0){50}}
\put(65,60){1}
\put(115,110){2}
\end{picture}}
\put(600,0){\begin{picture}(200,150)
\put(0,0){\framebox(50,50){}}
\put(50,0){\framebox(50,50){}}
\put(100,0){\framebox(50,50){}}
\put(50,50){\framebox(50,50){}}
\put(0,50){\framebox(50,50){}}
\put(0,0){\line(1,0){200}}
\put(0,100){\line(1,0){200}}
\put(0,0){\line(0,1){100}}
\put(200,0){\line(0,1){100}}
\put(150,0){\line(0,1){50}}
\put(150,50){\line(-1,0){50}}
\put(100,50){\line(0,1){50}}
\put(0,100){\framebox(150,50){}}
\put(165,10){1}
\put(115,110){2}
\end{picture}}
\end{picture}}

\bigskip\noindent
so that
$$
K_{443,\,221111111}(\zeta) = 2 K_{41,\,221}(\zeta) - K_{32,\,221}(\zeta)
= 2(\zeta^2 + \zeta^3) - (\zeta + \zeta^2) = 2\zeta^2 + 3\ .
$$
}
\end{example}

When $|\alpha|\le |\lambda_{(k)}|$, (\ref{KOSTROOT}) becomes simpler. For
if $|\alpha| < |\lambda_{(k)}|$ then $K_{\lambda,\,\nu}(\zeta) = 0$, and
otherwise the sum reduces to one single term
$$
K_{\lambda,\,\nu}(\zeta) = (-1)^{(k-1)|\mu|}
\epsilon_k(\lambda/\lambda_{(k)})\, K_{\lambda/\lambda_{(k)},\,\mu}^{(k)}\,
K_{\lambda_{(k)},\,\alpha}(\zeta) \ .
$$
In particular, for $\nu = (1^n)$, one recovers the expression of
$K_{\lambda,(1^n)}(\zeta)$ given by Morris and Sultana \cite{MS}.

%%%%%%%%%%%%%%%%%%%%%%%%%%%%%%%%%%%%%%%%%%%%%%%%%%%%%%%%%%%%%%%%%%%%%%%%
Finally, let us observe that the notion of $k$-sign of a partition can be
lifted
to a statistic on ribbon tableaux, which for technical reasons that
will become transparent later, takes values in $\N+{1\over 2}\N$, and
will be called {\it spin}.

Let $R$ be a $k$-ribbon, $h(R)$ its {\it heigth} and $w(R)$ its {\it width}.

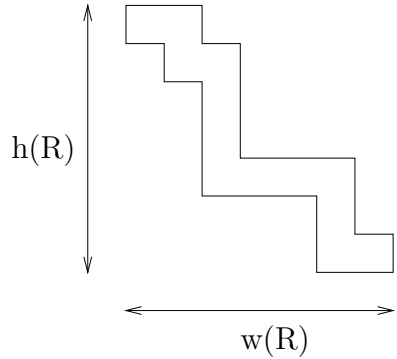
\begin{center}
\setlength{\unitlength}{0.01in}
\begin{picture}(200,195)(0,-10)
\path(42.000,172.000)(40.000,180.000)(38.000,172.000)
\path(40,180)(40,40)
\path(38.000,48.000)(40.000,40.000)(42.000,48.000)
\path(60,180)(100,180)(100,160)
        (120,160)(120,100)(180,100)
        (180,60)(200,60)(200,40)
        (160,40)(160,80)(100,80)
        (100,140)(80,140)(80,160)
        (60,160)(60,180)
\path(68.000,22.000)(60.000,20.000)(68.000,18.000)
\path(60,20)(200,20)
\path(192.000,18.000)(200.000,20.000)(192.000,22.000)
\put(120,0){\makebox(0,0)[lb]{\raisebox{0pt}[0pt][0pt]{\shortstack[l]{{\twlrm
w(R)}}}}}
\put(0,100){\makebox(0,0)[lb]{\raisebox{0pt}[0pt][0pt]{\shortstack[l]{{\twlrm
h(R)}}}}}
\end{picture}
\end{center}
The {\it spin} of $R$, denoted by $s(R)$, is defined as
\begin{equation}
s(R) ={h(R)-1\over 2}
\end{equation}
and the spin of a ribbon tableau $T$ is by definition the sum of the spins
of its ribbons. For example, the ribbon tableau

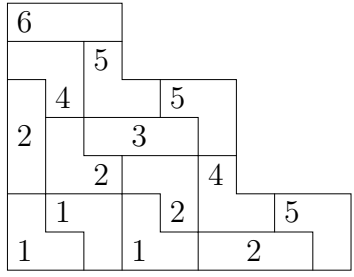
\begin{center}
\setlength{\unitlength}{0.01in}
\begingroup\makeatletter\ifx\SetFigFont\undefined
% extract first six characters in \fmtname
\def\x#1#2#3#4#5#6#7\relax{\def\x{#1#2#3#4#5#6}}%
\expandafter\x\fmtname xxxxxx\relax \def\y{splain}%
\ifx\x\y   % LaTeX or SliTeX?
\gdef\SetFigFont#1#2#3{%
  \ifnum #1<17\tiny\else \ifnum #1<20\small\else
  \ifnum #1<24\normalsize\else \ifnum #1<29\large\else
  \ifnum #1<34\Large\else \ifnum #1<41\LARGE\else
     \huge\fi\fi\fi\fi\fi\fi
  \csname #3\endcsname}%
\else
\gdef\SetFigFont#1#2#3{\begingroup
  \count@#1\relax \ifnum 25<\count@\count@25\fi
  \def\x{\endgroup\@setsize\SetFigFont{#2pt}}%
  \expandafter\x
    \csname \romannumeral\the\count@ pt\expandafter\endcsname
    \csname @\romannumeral\the\count@ pt\endcsname
  \csname #3\endcsname}%
\fi
\fi\endgroup
\begin{picture}(180,155)(0,-10)
\path(0,0)(180,0)(180,40)
        (120,40)(120,100)(60,100)
        (60,140)(0,140)(0,0)
\path(0,40)(20,40)(20,20)
        (40,20)(40,0)
\path(20,40)(60,40)(60,0)
\path(60,40)(80,40)(80,20)
        (100,20)(100,0)
\path(100,20)(160,20)(160,0)
\path(140,40)(140,20)
\path(0,100)(20,100)(20,40)
\path(20,80)(40,80)(40,60)
        (60,60)(60,40)
\path(60,60)(80,60)(100,60)(100,20)
\path(100,60)(120,60)
\path(0,120)(20,120)(40,120)(40,80)
\path(40,120)(60,120)
\path(40,80)(100,80)(100,60)
\path(80,100)(80,80)
\put(5,5){\makebox(0,0)[lb]{\smash{{{\SetFigFont{12}{14.4}{rm}1}}}}}
\put(5,65){\makebox(0,0)[lb]{\smash{{{\SetFigFont{12}{14.4}{rm}2}}}}}
\put(25,25){\makebox(0,0)[lb]{\smash{{{\SetFigFont{12}{14.4}{rm}1}}}}}
\put(45,45){\makebox(0,0)[lb]{\smash{{{\SetFigFont{12}{14.4}{rm}2}}}}}
\put(65,65){\makebox(0,0)[lb]{\smash{{{\SetFigFont{12}{14.4}{rm}3}}}}}
\put(85,85){\makebox(0,0)[lb]{\smash{{{\SetFigFont{12}{14.4}{rm}5}}}}}
\put(65,5){\makebox(0,0)[lb]{\smash{{{\SetFigFont{12}{14.4}{rm}1}}}}}
\put(85,25){\makebox(0,0)[lb]{\smash{{{\SetFigFont{12}{14.4}{rm}2}}}}}
\put(105,45){\makebox(0,0)[lb]{\smash{{{\SetFigFont{12}{14.4}{rm}4}}}}}
\put(125,5){\makebox(0,0)[lb]{\smash{{{\SetFigFont{12}{14.4}{rm}2}}}}}
\put(145,25){\makebox(0,0)[lb]{\smash{{{\SetFigFont{12}{14.4}{rm}5}}}}}
\put(25,85){\makebox(0,0)[lb]{\smash{{{\SetFigFont{12}{14.4}{rm}4}}}}}
\put(45,105){\makebox(0,0)[lb]{\smash{{{\SetFigFont{12}{14.4}{rm}5}}}}}
\put(5,125){\makebox(0,0)[lb]{\smash{{{\SetFigFont{12}{14.4}{rm}6}}}}}
\end{picture}
\end{center}
has a spin equal to $6$.

The $k$-sign of a partition $\lambda$ is thus equal to $(-1)^{2s(T)}$,
for any ribbon tableau $T$ of shape $\lambda$. For example, we can rewrite
the particular case $\nu=(0)$ of formula (\ref{plethrub}) as
\begin{equation}
\psi^k(h_\mu)=p_k\circ h_\mu = \sum_{T\in\tab_k(\,\cdot\, ,\mu)}(-1)^{2s(T)}s_T
\end{equation}
where $\tab_k(\,\cdot\, ,\mu)$ is the set of $k$-ribbon tableaux
of weight $\mu$, and
$s_T=s_\lambda$ if $\lambda$ is the shape of $T$.
We shall see in the next section that this formula leads to a simple
construction of a basis of highest weight vectors in the Fock space
representation of the quantum affine algebra $U_q(\slchap_n)$.

%%%%%%%%%%%%%%%%%%%%%%%%%%%%%%%%%%%%%%%%%%%%%%%%%%%%%%%%%%%%%%%%%%%%%%%%%%%
\section{The Fock representation of $U_q(\slchap_n)$}\label{S5}

The affine Lie algebra $\slchap_n=A_{n-1}^{(1)}$ has a natural
action on the space $\sym$ of symmetric functions, called the
bosonic Fock space representation. This representation is
equivalent to the infinite wedge, or fermionic Fock space representation,
and the isomorphism can be realized by means of vertex operators
(see  {\it e.g.} \cite{Kac}).

Let us recall briefly the fermionic version. Let $V$ be the vector
space $\C^{(\Z)}$ and $(u_i)_{i\in\Z}$ be its canonical basis.
The fermionic Fock space $\bigwedge^\infty V$ is defined as the vector
space spanned by the infinite exterior products
$u_{i_1}\wedge u_{i_2}\wedge\cdots\wedge u_{i_n}\wedge\cdots $ satisfying
$i_1> i_2 >i_3>\ldots > i_n>\ldots $ and $i_k-i_{k+1}=1$ for $k>>0$.
The wedge product is as usual alternating, and linear in each
of its factors. We denote by $\F$ the subspace of $\bigwedge^\infty V$
spanned by the elements such that $i_k=-k+1$ for $k>>0$
(usually this subspace is denoted by $\F^{(0)}$, but as we
shall not need the other sectors $\F^{(m)}$, we drop the superscript).

The Lie algebra $\gl_\infty$ of $\Z\times \Z$ complex matrices
$A=(a_{ij})$ with finitely many nonzero entries acts on $\bigwedge^\infty V$
by
$$
A\cdot u_{i_1}\wedge u_{i_2}\wedge\cdots =
(Au_{i_1})\wedge u_{i_2}\wedge\cdots
+u_{i_1}\wedge (Au_{i_2})\wedge\cdots + \cdots
$$
the sum having only a finite number of nonzero terms.

Let $E_{ij}$ be the infinite matrix $(E_{ij})_{rs}=\delta_{ir}\delta_{js}$.
The subalgebra $\Sl_\infty$ of $\gl_\infty$ constituted by the matrices
with zero trace has for Chevalley generators $e_i^\infty =E_{i,i+1}$,
$f_i^\infty=E_{i+1,i}$ and $h_i=E_{i,i}-E_{i+1,i+1}$, $i\in\Z$.
The infinite sums
$$
e_i = \sum_{j\equiv i \mod n} e_j^\infty\ ,  \qquad
f_i = \sum_{j\equiv i \mod n} f_j^\infty\ ,  \qquad
h_i = \sum_{j\equiv i \mod n} h_j^\infty\
$$
do not belong to $\gl_\infty$, but they have a well-defined action
on $\bigwedge^\infty V$, and it can be checked that they generate a
representation of $\widehat{\Sl}'_n$ with central charge $c=1$.
The remaining generator $D$ of $\slchap_n =\widehat{\Sl}'_n\oplus\C D$
can be implemented by
$$
D:=-\sum_{i\in\Z} \left[ {i\over n}\right] (E_{ii}-\theta(-i-1)) \ ,
$$
where $\theta(x) = 1$ for $x\ge 0$ and $\theta(x)=0$ otherwise
(\cf \cite{DJKMO}).

The basis vectors of $\F$ can be labelled by
partitions, by setting
$$
|\lambda\> = u_{\lambda_1}\wedge u_{\lambda_2-1}\wedge
u_{\lambda_3-2}\wedge\cdots \ .
$$
With this indexation, the action of the Chevalley generators of
$\slchap_n$ can be described as follows \cite{DJKMO}.

To each node $(i,j)$ of a Young diagram, one can associate its
{\it residue} $\rho_{i,j} = j-i \mod n \in\{0,\ldots,n-1\}$.
Then,
\begin{equation}\label{rIND}
e_r|\lambda\> =\sum |\mu\> \ , \qquad f_r|\lambda\>=\sum |\nu\>
\end{equation}
where $\mu$ ({\it resp.} $\nu$) runs  over all diagrams obtained
from $\lambda$ by removing ({\it resp. } adding) a node of residue $r$.

In this picture, one can observe that $e_r$ and $f_r$ are exactly
the $r$-restricting and $r$-inducing operators introduced by
G. de B. Robinson in the context of the modular representation
theory of the symmetric group \cite{Ro}.

The natural way of interpreting the basis vectors $|\lambda\>$ as symmetric
functions is to put $|\lambda\>=s_\lambda$. This is imposed by the
boson-fermion correspondence, and it is also compatible with the
modular representation interpretation.

In this realization, it can be shown that the image $U(\slchap_n)|0\>$
of the constant $|0\>=s_0=1$, which is the  basic representation
$M(\Lambda_0)$, is equal to the subalgebra
$$
{\cal T}^{(n)}=\C[p_i\, |\, i\not\equiv 0 \mod n]
$$
generated by the power-sums $p_i$ such that $n\not | \, i$.

The bosonic operators
$$
b_k\ :\ f \longmapsto  D_{p_{kn}}f =kn{\partial\over \partial p_{kn}} f
\quad {\rm and} \quad b_{-k}\ :\ f\longmapsto p_{kn}f \qquad (k\ge 1)
$$
commute with the action of $\slchap_n$. They generate a Heisenberg algebra
$\H$, and the irreducible $\H$-module $U(\H)|0\>$ is exactly the
space $\T^{(n)}$ of highest weight vectors of the Fock space, viewed as an
$\slchap_n$-module.
Thus, these highest weight vectors are exactly the plethysms
$\psi^n(f), f\in\sym$. Natural bases of $\T^{(n)}$ are therefore
$\psi^n(p_\mu)=p_{n\mu}$, $\psi^n(s_\mu)$ or $\psi^n(h_\mu)$.
We know from Section \ref{S4} that this last one admits a
simple combinatorial description in terms of ribbon tableaux:
\begin{equation}\label{pnh}
\psi^n(h_\mu) = \sum_{\lambda}
\epsilon_n(\lambda/\mu)K_{\lambda\mu}^{(n)}s_\lambda
=\sum_{T\in\tab_n(\,\cdot\, ,\mu)} (-1)^{2s(T)}s_T \ .
\end{equation}

This  formula is especially meaningful in the quantized version,
that we shall now describe.

We first recall the definition of $U_q(\widehat{\Sl}_n)$
(\cf \cite{MM} and references therein). Let $\h$ be a
$(n+1)$-dimensional vector space over $\Q$ with basis
$\{h_0,h_1,\ldots ,h_{n-1},D\}$. We denote by
$\{\Lambda_0,\Lambda_1,\ldots ,\Lambda_{n-1},\delta\}$ the dual basis of
$\h^*$, that is,
$$
\<\Lambda_i,h_j\> = \delta_{ij}, \quad
\<\Lambda_i,D\> = 0,\quad
\<\delta, h_i\> = 0, \quad
\<\delta, D\> = 1 ,
$$
and we set
$\alpha_i = 2\Lambda_i - \Lambda_{i-1} - \Lambda_{i+1} + \delta_{i0} \,\delta$
for $i=0,1,\ldots ,n-1$. In these formulas it is understood that
$\Lambda_n = \Lambda_0$ and $\Lambda_{-1} = \Lambda_{n-1}$.
The $n\times n$ matrix $[\<\alpha_i,h_j\>]$ is the generalized Cartan
matrix associated to $\widehat{\Sl}_n$. The weight lattice is
$P=(\bigoplus_{i=0}^{n-1} \Z\Lambda_i) \bigoplus \Z\delta$, its dual
is $P^\vee = (\bigoplus_{i=0}^{n-1} \Z\,h_i) \bigoplus \Z\,D$,
and the root lattice is
$Q=\bigoplus_{i=0}^{n-1} \Z\alpha_i$.
One defines $U_q(\widehat{\Sl}_n)$ as the associative algebra with 1 over
$\Q(q)$
generated by the symbols $e_i,\ f_i,\ 0\le i \le n-1,$ and $q^h,\ h \in
P^\vee$,
subject to the relations
$$
q^h\,q^{h'}  =   q^{h+h'},  \quad  q^0 = 1,
$$
$$
q^h e_j q^{-h} = q^{\<\alpha_j,h\>} e_j ,
$$
$$
q^h f_j q^{-h} = q^{-\<\alpha_j,h\>} f_j ,
$$
$$
[e_i,f_j] = \delta_{ij} {q^{h_i}- q^{-h_i} \over q - q^{-1}} ,
$$
$$
\sum_{k=0}^{1-\<\alpha_i,h_j\>} (-1)^k
\left[ \begin{array}{c}
  1-\<\alpha_i,h_j\> \\ k
  \end{array}
\right]
e_i^{1-\<\alpha_i,h_j\> -k} e_j e_i^k = 0  \quad (i\not = j) ,
$$
$$
\sum_{k=0}^{1-\<\alpha_i,h_j\>} (-1)^k
\left[ \begin{array}{c}
  1-\<\alpha_i,h_j\> \\ k
	\end{array}
	\right]
	f_i^{1-\<\alpha_i,h_j\> -k} f_j f_i^k = 0  \quad (i\not = j) .
$$
Here  the  $q$-integers,
$q$-factorials and $q$-binomial coefficients are the symmetric ones:
$$
[k] = {q^k- q^{-k} \over q - q^{-1}} , \quad
[k]! = [k]\,[k-1]\,\cdots [1] , \quad
\left[
\begin{array}{c}
m \\ k
\end{array}
\right]
= {[m]!\over [m-k]!\,[k]!} \,.
$$

We now recall some definitions relative to $U_q(\widehat{\Sl}_n)$-modules.
Let $M$ be a $U_q(\widehat{\Sl}_n)$-module and $\Lambda \in P$ a weight.
The subspace
$$M_\Lambda = \{ v\in M \ | \ q^h\,v = q^{\<\Lambda,h\>} \,v, \ h\in P^\vee
\}$$
is called the weight space of weight $\Lambda$ of $M$ and its elements
are called the weight vectors of weight $\Lambda$.
The module $M$ is said to be integrable if
\begin{quote}
(i)\quad $M=\bigoplus_{\Lambda \in P} M_\Lambda ,$

(ii)\quad ${\rm dim}\,M_\Lambda < \infty$ for $\Lambda \in P$,

(iii)\quad for $i= 0,1,\ldots ,n-1,$ $M$ decomposes into a direct sum of finite
dimensional $U_i$-modules, where $U_i$ denotes the subalgebra
of $U_q(\widehat{\Sl}_n)$ generated by $e_i,\ f_i,\ q^{h_i},\ q^{-h_i}.$
\end{quote}
A highest weight vector $v\in M$ is a vector annihilated by all raising
operators $e_i$. The module $M$ is said to be a highest weight
module if there exists a highest weight vector $v$ such that
$M = U_q(\widehat{\Sl}_n)\, v$. The weight of $v$ is called the highest
weight of $M$.

By the representation theory of $U_q(\widehat{\Sl}_n)$, there exists for each
dominant integral weight $\Lambda$ ({\it i.e.} $\<\Lambda , h_i\> \in \Z_+$
for $i= 0,1,\ldots ,n-1$) a unique integrable highest weight module
$M(\Lambda)$ with highest weight $\Lambda$.

A $q$-analog of the Fock representation of $\slchap_n$ can be realized
in the $\Q(q)$-vector space $\F$ spanned by all partitions:
$$ \F = \bigoplus_{\lambda \in \P}  \Q(q) \,|\lambda\> \,$$
the action being defined in combinatorial terms.

Let us say that a point $(a,b)$ of $\Z_+\times \Z_+$ is an indent
$i$-node of a Young diagram $\lambda$ if a box of residue $i=a-b\mod n$
can be added to $\lambda$ at position $(a,b)$, in such a way that the
new diagram still corresponds to a partition.
Similarly, a node of $\lambda$ of residue $i$
which can be removed will be called a removable $i$-node.

Let $i\in \{0,1,\ldots ,n-1\}$ and let $\lambda$, $\nu$ be two partitions
such that $\nu$ is obtained from $\lambda$ by filling an indent $i$-node
$\gamma$
We set:
\begin{quote}
$N_i(\lambda) = \sharp \{$ indent $i$-nodes of $\lambda$
$\} - \sharp \{$ removable $i$-nodes of $\lambda$ $\}$,

$N_i^l(\lambda,\nu) = \sharp \{$ indent $i$-nodes of $\lambda$ situated to
the {\it left} of $\gamma$ (not counting $\gamma$) $\}$
$- \sharp \{$ removable $i$-nodes of $\lambda$ situated
to the {\it left} of $\gamma$ $\}$,

$N_i^r(\lambda,\nu) = \sharp \{$ indent $i$-nodes of $\lambda$ situated to
the {\it right} of $\gamma$ (not counting $\gamma)~\} - \sharp \{$
removable $i$-nodes of $\lambda$ situated to the {\it right} of $\gamma$ $\}$,

$N^0(\lambda) = \sharp \{$ 0-nodes of $\lambda$ $ \}$.
\end{quote}
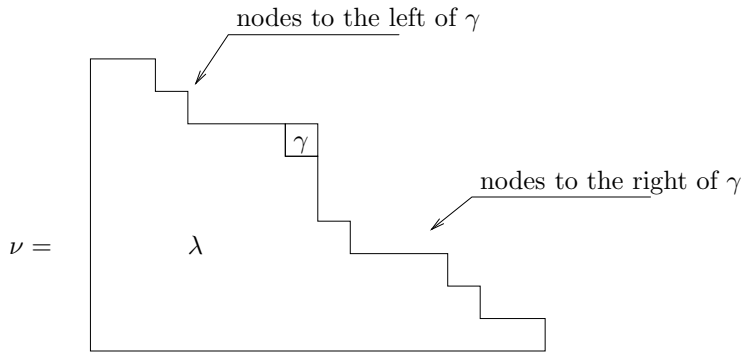
\begin{center}

\setlength{\unitlength}{0.0085in}
\begingroup\makeatletter\ifx\SetFigFont\undefined
% extract first six characters in \fmtname
\def\x#1#2#3#4#5#6#7\relax{\def\x{#1#2#3#4#5#6}}%
\expandafter\x\fmtname xxxxxx\relax \def\y{splain}%
\ifx\x\y   % LaTeX or SliTeX?
\gdef\SetFigFont#1#2#3{%
  \ifnum #1<17\tiny\else \ifnum #1<20\small\else
  \ifnum #1<24\normalsize\else \ifnum #1<29\large\else
  \ifnum #1<34\Large\else \ifnum #1<41\LARGE\else
     \huge\fi\fi\fi\fi\fi\fi
  \csname #3\endcsname}%
\else
\gdef\SetFigFont#1#2#3{\begingroup
  \count@#1\relax \ifnum 25<\count@\count@25\fi
  \def\x{\endgroup\@setsize\SetFigFont{#2pt}}%
  \expandafter\x
    \csname \romannumeral\the\count@ pt\expandafter\endcsname
    \csname @\romannumeral\the\count@ pt\endcsname
  \csname #3\endcsname}%
\fi
\fi\endgroup
\begin{picture}(399,230)(0,-10)
\path(50,180)(50,0)(330,0)
	(330,20)(290,20)(290,40)
	(270,40)(270,60)(210,60)
	(210,80)(190,80)(190,100)
	(190,120)(170,120)(170,140)
	(110,140)(110,160)(90,160)
	(90,180)(50,180)
\path(190,120)(190,140)(170,140)
	(170,120)(190,120)
\path(240,195)(135,195)(115,165)
\path(117.774,172.766)(115.000,165.000)(121.102,170.547)
\path(395,95)(285,95)(260,75)
\path(264.998,81.559)(260.000,75.000)(267.496,78.436)
\put(0,60){\makebox(0,0)[lb]{\smash{{{\SetFigFont{10}{14.4}{rm}$\nu =$}}}}}
\put(110,60){\makebox(0,0)[lb]{\smash{{{\SetFigFont{10}{14.4}{rm}$\lambda$}}}}}
\put(175,125){\makebox(0,0)[lb]{\smash{{{\SetFigFont{10}{14.4}{rm}$\gamma$}}}}}
\put(140,200){\makebox(0,0)[lb]{\smash{{{\SetFigFont{10}{14.4}{rm}nodes to the
left of $\gamma$}}}}}
\put(290,100){\makebox(0,0)[lb]{\smash{{{\SetFigFont{10}{14.4}{rm}nodes to the
right of $\gamma$}}}}}
\end{picture}
\end{center}
The following result is due to Hayashi \cite{Hay}, and the formulation
that we use has been given by Misra and Miwa \cite{MM}
(with a slight change in the conventions, that is conjugation
of partitions and $q\rightarrow 1/q$).

\begin{theorem}\label{HAYASHI}
The algebra $U_q(\widehat{\Sl}_n)$ acts on $\F$ by
\begin{quote}
$q^{h_i} \,|\lambda\> = q^{N_i(\lambda)}\,|\lambda\>\,,$

$q^D \, |\lambda\> = q^{-N^0(\lambda)} \, |\lambda\>\,,$

$f_i |\lambda\> = \sum_\nu q^{N_i^r(\lambda,\nu)} \, |\nu\> \,,$
sum over all partitions $\nu$ such that $\nu/\lambda$ is a $i$-node,

$e_i |\nu\> = \sum_\lambda q^{-N_i^l(\lambda,\nu)} \, |\lambda\> \,,$
sum over all partitions $\lambda$ such that $\nu/\lambda$ is a $i$-node.
\end{quote}
\end{theorem}
It is easy to see that $\F$ is an integrable $U_q(\widehat{\Sl}_n)$-module.
It is not irreducible.
Actually it decomposes as
$$\F \cong \bigoplus_{k\ge 0} M(\Lambda_0 - k\delta)^{\oplus p(k)} \,.$$
Obviously, the empty partition $|0\>$ is a highest weight vector
of weight $\Lambda_0$. The submodule $U_q(\widehat{\Sl}_n) \, |0\>$ is
isomorphic to $M(\Lambda_0)$, also called the {\it basic representation} of
$U_q(\widehat{\Sl}_n)$.
Again, one can identify $\F$ with $\sym$ (with coefficients in $\Q(q)$ and
interpret $|\lambda\>$ as $s_\lambda$.
Then, a natural $q$-analog of (\ref{pnh}) gives a basis
of highest weight vectors for $U_q(\slchap_n)$ in $\F$:

\begin{proposition}
Define a linear operator  $\psi^n_q$ on $\sym$ by
$$
\psi_q^n(h_\mu)=\sum_{T\in\tab_n(\,\cdot\, ,\mu)}
(-q)^{-2s(T)} s_T \ .
$$
Then, its  image $\psi_q^n(\sym)$ is the space ${\cal S}^{(n)}_q$
of highest weight vectors of $U_q(\slchap_n)$ in $\sym$.
\end{proposition}

\begin{example}{\rm The plethysm $\psi^2(h_{21})$ is given by the
following domino tableaux

\bigskip
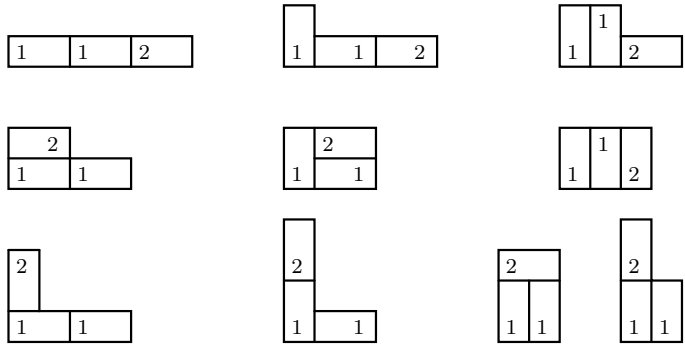
\begin{center}
%avec xfig du litp:
%\setlength{\unitlength}{0.00825in}

%\setlength{\unitlength}{0.00083333in}
%pour monge
\setlength{\unitlength}{0.00053333in}
\begingroup\makeatletter\ifx\SetFigFont\undefined
% extract first six characters in \fmtname
\def\x#1#2#3#4#5#6#7\relax{\def\x{#1#2#3#4#5#6}}%
\expandafter\x\fmtname xxxxxx\relax \def\y{splain}%
\ifx\x\y   % LaTeX or SliTeX?
\gdef\SetFigFont#1#2#3{%
  \ifnum #1<17\tiny\else \ifnum #1<20\small\else
  \ifnum #1<24\normalsize\else \ifnum #1<29\large\else
  \ifnum #1<34\Large\else \ifnum #1<41\LARGE\else
     \huge\fi\fi\fi\fi\fi\fi
  \csname #3\endcsname}%
\else
\gdef\SetFigFont#1#2#3{\begingroup
  \count@#1\relax \ifnum 25<\count@\count@25\fi
  \def\x{\endgroup\@setsize\SetFigFont{#2pt}}%
  \expandafter\x
    \csname \romannumeral\the\count@ pt\expandafter\endcsname
    \csname @\romannumeral\the\count@ pt\endcsname
  \csname #3\endcsname}%
\fi
\fi\endgroup
\begin{picture}(6624,3339)(0,-10)
\thicklines
\path(612,3012)(612,2712)
\path(1212,3012)(1212,2712)
\path(12,2112)(12,1512)(612,1512)
	(1212,1512)(1212,1812)(612,1812)
	(612,2112)(12,2112)
\path(12,1812)(612,1812)(612,1512)
\path(12,912)(312,912)
\path(312,912)(312,312)(1212,312)
	(1212,12)(12,12)(12,912)
\path(12,312)(312,312)
\path(612,312)(612,12)
\path(2712,3312)(2712,2712)(4212,2712)
	(4212,3012)(3012,3012)(3012,3312)(2712,3312)
\path(3012,3012)(3012,2712)
\path(3612,3012)(3612,2712)
\path(2712,2112)(2712,1512)(3612,1512)
	(3612,2112)(2712,2112)
\path(3012,2112)(3012,1512)
\path(3012,1812)(3612,1812)
\path(2712,1212)(2712,12)(3612,12)
	(3612,312)(3012,312)(3012,1212)(2712,1212)
\path(2712,612)(3012,612)
\path(3012,312)(3012,12)
\path(5412,3312)(5412,2712)(6612,2712)
	(6612,3012)(6012,3012)(6012,3312)(5412,3312)
\path(5712,3312)(5712,2712)
\path(6012,3012)(6012,2712)
\path(5412,2112)(5412,1512)(6312,1512)
	(6312,2112)(5412,2112)
\path(5712,2112)(5712,1512)
\path(6012,2112)(6012,1512)
\path(4812,912)(4812,12)(5412,12)
	(5412,912)(4812,912)
\path(6012,1212)(6012,12)(6612,12)
	(6612,612)(6312,612)(6312,1212)(6012,1212)
\path(6012,612)(6312,612)(6312,12)
\path(4812,612)(5412,612)
\path(5112,612)(5112,12)
\path(12,3012)(12,2712)(612,2712)
	(1212,2712)(1812,2712)(1812,3012)(12,3012)
\put(87,2787){\makebox(0,0)[lb]{\smash{{{\SetFigFont{8}{9.6}{rm}1}}}}}
\put(3987,2787){\makebox(0,0)[lb]{\smash{{{\SetFigFont{8}{9.6}{rm}2}}}}}
\put(687,2787){\makebox(0,0)[lb]{\smash{{{\SetFigFont{8}{9.6}{rm}1}}}}}
\put(1287,2787){\makebox(0,0)[lb]{\smash{{{\SetFigFont{8}{9.6}{rm}2}}}}}
\put(87,1587){\makebox(0,0)[lb]{\smash{{{\SetFigFont{8}{9.6}{rm}1}}}}}
\put(387,1887){\makebox(0,0)[lb]{\smash{{{\SetFigFont{8}{9.6}{rm}2}}}}}
\put(687,1587){\makebox(0,0)[lb]{\smash{{{\SetFigFont{8}{9.6}{rm}1}}}}}
\put(87,87){\makebox(0,0)[lb]{\smash{{{\SetFigFont{8}{9.6}{rm}1}}}}}
\put(687,87){\makebox(0,0)[lb]{\smash{{{\SetFigFont{8}{9.6}{rm}1}}}}}
\put(87,687){\makebox(0,0)[lb]{\smash{{{\SetFigFont{8}{9.6}{rm}2}}}}}
\put(2787,2787){\makebox(0,0)[lb]{\smash{{{\SetFigFont{8}{9.6}{rm}1}}}}}
\put(2787,1587){\makebox(0,0)[lb]{\smash{{{\SetFigFont{8}{9.6}{rm}1}}}}}
\put(3387,1587){\makebox(0,0)[lb]{\smash{{{\SetFigFont{8}{9.6}{rm}1}}}}}
\put(3087,1887){\makebox(0,0)[lb]{\smash{{{\SetFigFont{8}{9.6}{rm}2}}}}}
\put(2787,87){\makebox(0,0)[lb]{\smash{{{\SetFigFont{8}{9.6}{rm}1}}}}}
\put(3387,87){\makebox(0,0)[lb]{\smash{{{\SetFigFont{8}{9.6}{rm}1}}}}}
\put(2787,687){\makebox(0,0)[lb]{\smash{{{\SetFigFont{8}{9.6}{rm}2}}}}}
\put(4887,87){\makebox(0,0)[lb]{\smash{{{\SetFigFont{8}{9.6}{rm}1}}}}}
\put(5187,87){\makebox(0,0)[lb]{\smash{{{\SetFigFont{8}{9.6}{rm}1}}}}}
\put(4887,687){\makebox(0,0)[lb]{\smash{{{\SetFigFont{8}{9.6}{rm}2}}}}}
\put(6087,87){\makebox(0,0)[lb]{\smash{{{\SetFigFont{8}{9.6}{rm}1}}}}}
\put(6387,87){\makebox(0,0)[lb]{\smash{{{\SetFigFont{8}{9.6}{rm}1}}}}}
\put(6087,687){\makebox(0,0)[lb]{\smash{{{\SetFigFont{8}{9.6}{rm}2}}}}}
\put(5487,1587){\makebox(0,0)[lb]{\smash{{{\SetFigFont{8}{9.6}{rm}1}}}}}
\put(5787,1887){\makebox(0,0)[lb]{\smash{{{\SetFigFont{8}{9.6}{rm}1}}}}}
\put(6087,1587){\makebox(0,0)[lb]{\smash{{{\SetFigFont{8}{9.6}{rm}2}}}}}
\put(5487,2787){\makebox(0,0)[lb]{\smash{{{\SetFigFont{8}{9.6}{rm}1}}}}}
\put(5787,3087){\makebox(0,0)[lb]{\smash{{{\SetFigFont{8}{9.6}{rm}1}}}}}
\put(6087,2787){\makebox(0,0)[lb]{\smash{{{\SetFigFont{8}{9.6}{rm}2}}}}}
\put(3387,2787){\makebox(0,0)[lb]{\smash{{{\SetFigFont{8}{9.6}{rm}1}}}}}
\end{picture}
\end{center}

\medskip\noindent
and the corresponding highest weight vector of $U_q(\slchap_2)$ is
$$
\psi^2_q(h_{21})=
s_{6}-q^{-1}s_{51}+(1+q^{-2})s_{42}-q^{-1}s_{411}$$
$$
-(q^{-1}+q^{-3})s_{33}+q^{-2}s_{3111}+q^{-2}s_{222}-q^{-3}s_{2211}
$$

}
\end{example}

The proposition is a consequence of  the following more precise statement.

\begin{theorem}\label{UV1}
Let $U_k$, $V_k$ ($k\ge 1$) be the linear operators defined by
$$
V_k s_\lambda =
\sum_\mu \left(
\sum_{T\in\tab_n(\mu/\lambda,(k))} (-q)^{-2s(T)} \right) s_\mu \ ,
$$
$$
U_k s_\lambda =
\sum_\nu \left(
\sum_{T\in\tab_n(\lambda/\nu,(k))} (-q)^{-2s(T)} \right) s_\nu \ ,
$$
so that $V_k$ is a $q$-analog of $f\mapsto \psi^n(h_k)f$, and $U_k$ is
its adjoint. Then, $U_k$ and $V_k$ commute with the action of
$U_q(\slchap_n)$. In particular, each
$\psi^n_q(h_\mu)=V_{\mu_r}\cdots V_{\mu_1}|0\>$
is a highest weight vector.
\end{theorem}

This result can be obtained by a direct verification, using
formula (\ref{plethrub}). However, a more illuminating approach comes
from comparison with a recent construction of Stern \cite{Ste} and
Kashiwara, Miwa and Stern \cite{KMS}. These authors construct the
$q$-analog of the Fock representation by means of a $q$-deformation
of the wedge product, defined in terms of an action of the affine
Hecke algebra $\widehat{H}_N(q^{-2})$ on a tensor product
$V(z)^{\otimes n}$ of evaluation modules.

Here, $V(z)$ is $\C^{(\Z)}$ realized as
$$
\left( \bigoplus_{i=1}^n \C\, v_i \right) \otimes \C[z,z^{-1}] \ ,
$$
where $z^iv_j$ is identified with $u_{j-ni}$, endowed with an
appropriate $q$-analog of the action of $\slchap_n$ on $V$.

Writing $z^{r_1}v_{m_1}\otimes\cdots\otimes z^{r_N}v_{m_N}$
as
$v_{{\bf m}}z^{{\bf r}}=v_{m_1}\otimes\cdots\otimes v_{m_N}\cdot
z_1^{r_1}\cdots z_N^{r_N}$,
the right action of $\widehat{H}_N(q^{-2})$ on $V(z)^{\otimes N}$ is
described by the following formulas \cite{GRV,Ste,KMS}:
\begin{quote}
$y_i$ acts as $z_i^{-1}$ and
\end{quote}
\begin{equation}\label{ACTHECKE}
(v_{{\bf m}}\cdot z^{{\bf r}})T_i =
\left\{\matrix{ -q^{-1} v_{{\bf m}\sigma_i}\cdot \sigma_i(z^{{\bf r}})
		+(q^{-2}-1)v_{{\bf m}}\cdot \partial_i(z_i z^{{\bf r}})\ {\rm if}\
		m_i<m_{i+1}  \cr
		- v_{{\bf m}}\cdot \sigma_i(z^{{\bf r}})
		+(q^{-2}-1)v_{{\bf m}}\cdot z_i\partial_i(z^{{\bf r}})\ {\rm if}\
		m_i=m_{i+1}   \cr
		-q^{-1}v_{{\bf m}\sigma_i}\cdot \sigma_i(z^{{\bf r}})
		+(q^{-2}-1)v_{{\bf m}}\cdot z_i\partial_i(z^{{\bf r}})\ {\rm if}\
		m_i>m_{i+1} \cr
		}\right.
\end{equation}
where ${\bf m}\sigma_i=(m_{\sigma(1)},\ldots,m_{\sigma(N)})$ and
$\partial_i$ is the divided difference operator
$f(z)\mapsto (f-\sigma_i(f))/(z_i-z_{i+1})$. This action can be regarded
as a generalization of the one given in \cite{DKLLST}, which
would correspond to the degenerate case $n=1$.

The important point is that this action commutes with $U_q(\slchap_n)$.
Let $A^{(N)}=\sum_{\sigma\in\S_N} T_\sigma$. This is a $q$-analog
of the total antisymmetrizer of $\S_N$, since  signs have
been incorporated in formulas (\ref{ACTHECKE}) in such a way that
$T_i$ acts as a $q$-analog of $-\sigma_i$. Kashiwara, Miwa and Stern
define then the $q$-exterior powers by
$\bigwed^N V(z)=V(z)/\ker A^{(N)}$, and denote by
$u_{i_1}\wed\cdots\wed u_{i_N}$ the image of
$u_{i_1}\otimes \cdots\otimes u_{i_N}$ in the quotient.
A basis of $\bigwed^N V(z)$ is formed by the normally ordered products
$u_{i_1}\wed\cdots\wed u_{i_N}$, where $i_1>i_2>\ldots >i_N$,
and any $q$-wedge product can be expressed on this basis, by means
of the following relations iteratively applied to consecutive factors.
Suppose that $\ell< m$ and that $\ell-m \mod n = i$. Set $t=q^{-1}$.
Then,
\begin{quote}
--- if $i=0$ then $u_\ell\wed u_m = -u_m \wed u_\ell$ \\
--- otherwise, $u_\ell\wed u_m = -t u_m\wed u_\ell
+(t^2-1)(u_{m-i}\wed u_{\ell+i}-t u_{m-n}\wed u_{m+n}
+t^2 u_{m-n-i}\wed u_{l+n+i} +\cdots)$
\end{quote}
where the only terms to be taken into account in this last expression
are the normally ordered ones.

There is then a well-defined action of $U_q(\slchap_n)$
on the ``thermodynamic limit''
$\F_q=\bigwed^\infty V(z)=\lim_{N\rightarrow \infty}\bigwed^N V(z)$,
which provides another realization of the $q$-Fock representation.

The affine Hecke algebra does not act anymore on $\bigwed^N V(z)$, but
its center does. This center is generated by the power sums
$$
p_k(Y) =\sum_{i=1}^{N} y_i^k \qquad k=\pm 1,\pm 2,\ldots
$$
and at the thermodynamic limit, the operators
$$
B_k =\sum_{i=1}^\infty y_i^{-k}
$$
are shown in \cite{KMS} to generate an action of a Heisenberg algebra
on $\F_q$, with
\begin{equation}
[B_k,B_\ell] = k {1-q^{2nk}\over 1-q^{2k}}\cdot \delta_{k,-\ell} \ .
\end{equation}

If one interprets the infinite $q$-wedges as Schur functions, by
the same rule as in the classical case, one sees that $B_{-k}$  ($k\ge 1$) is
a $q$-analogue of the multiplication operator $f\mapsto p_{nk}f$,
and that $B_k$ corresponds to its ajoint $D_{p_{nk}}$.

To connect this construction with the preceding one, take for generators
of the center of the affine Hecke algebra the elementary symmetric
functions in the $y_i$ and $y_i^{-1}$ instead of the power sums, and
define operators on $\F_q$ by
$$
\tilde U_k = e_k(y_1,y_2,\ldots ),\qquad
\tilde V_k = e_k(y_1^{-1},y_2^{-1},\ldots ) \ .
$$
These operators commute with  $U_q(\slchap_n)$, and their
action can be described in terms of ribbon tableaux:

\begin{lemma}
\begin{equation}
\tilde U_k|\lambda\>=
\sum_\nu \left( \sum_{T\in \tab_n(\lambda'/\nu',(k))}(-q)^{-2s(T)}
 \right)|\nu\>
\end{equation}
\begin{equation}
\tilde V_k|\lambda\> =
\sum_\mu \left( \sum_{T\in\tab_n(\mu'/\lambda',(k))}(-q)^{-2s(T)}
\right) |\mu\> \ .
\end{equation}
\end{lemma}

\Proof It is sufficient to work with $\bigwed^NV(z)$ for $N$ sufficiently
large. Then,
$$
\tilde V_k u_{i_1}\wed\cdots\wed u_{i_N} = \sum_J
u_{i_1+j_1}\wed\cdots\wed u_{i_N+j_N}
$$
where $J$ runs through the distinct permutations of the integer
vector $(0^{N-k} n^k)$. The only reorderings needed to express a
term of this sum in standard form are due to the appearance of factors
of the form $u_i\wed u_{j+n}$ with $j+n>i>j$. In this case,
$u_i\wed u_{j+n}= -t u_{j+n}\wed u_i$ since the other terms
$(t^2-1)(u_{j+n-a}\wed u_{i+a}-tu_{j}\wed u_{i+n}+\cdots)$ vanish,
the residue $a=j+n-i\mod n$ being actually equal to $j+n-i$.
The first case of the straightening rule is never encountered
because $j+n-i\equiv 0\mod n$ would imply $i-j=bn$ with $b>0$,
so that $j+n \not > i$.

Thus,
$$
u_{i_1+j_1}\wed\cdots\wed u_{i_N+j_N}
=(-t)^{\ell(\sigma)} u_{i_{\sigma(1)}+j_{\sigma(1)}}\wed\cdots\wed
 u_{i_{\sigma(N)}+j_{\sigma(N)}}
$$
where $\sigma$ is the shortest permutation such that the result
is normally ordered. In view of the remarks in Section \ref{S4},
this gives the result for $\tilde V_k$. The argument for
$\tilde U_k$ is similar.

\begin{corollary}\label{UV}
The operators $U_k,V_k$ of Theorem \ref{UV1} act on $\F_q$
as $h_k(y_1,y_2,\ldots)$ and $h_k(y_1^{-1},y_2^{-1},\ldots)$ respectively.
In particular, $[U_i,U_j]=[V_i,V_j]=0$.
\end{corollary}

%%%%%%%%%%%%%%%%%%%%%%%%%%%%%%%%%%%%%%%%%%%%%%%%%%%%%%%%%%%%%%%%%%%%

\section{$H$-functions}\label{S6}

Let $\lambda$ be a partition without $k$-core, and with $k$-quotient
$(\lambda^0,\ldots,\lambda^{k-1})$. For a ribbon tableau $T$ of
weight $\mu$, let $x^T=x_1^{\mu_1}x_2^{\mu_2}\cdots x_r^{\mu_r}$.
Then, the  correspondence between $k$-ribbon tableaux and
$k$-tuples of ordinary tableaux
shows that the generating function
\begin{equation}
\G^{(k)}_\lambda
= \sum_{T\in \tab_k(\lambda,\,\cdot\,)}x^T
= \prod_{i=0}^{k-1}\ \sum_{\t_i\in\tab(\lambda^i,\,\cdot\,)}x^{\t_i}
=\prod_{i=0}^{k-1} s_{\lambda^i}
\end{equation}
is a product of Schur functions. Introducing in this equation an appropriate
statistic
on ribbon tableaux, one can therefore obtain $q$-analogues of products of Schur
functions.
The statistic called {\it cospin},  described below, leads to $q$-analogues
with interesting properties.

For a partition $\lambda$ without $k$-core, let
\begin{equation}
s^*_k(\lambda)=\max \{s(T)\,|\, T\in \tab_k(\lambda,\,\cdots\,)\} \ .
\end{equation}
The {\it cospin} $\cs(T)$ of a $k$-ribbon tableau $T$ of shape $\lambda$
is then
\begin{equation}
\cs(T)=s^*_k(\lambda)-s(T) \ .
\end{equation}
Although $s(T)$ can be a half-integer, it is easily seen that $\cs(T)$
is always an integer. Also, there is one important case where $s(T)$
is an integer. This is when the shape $\lambda$ of $T$ is of the
form $k\mu=(k\mu_1,k\mu_2,\ldots,k\mu_r)$. In this case, the partitions
constituting the $k$-quotient of $\lambda$ are formed by parts of $\mu$,
grouped according to the class modulo $k$ of their indices. More
precisely, $\lambda^i=\{\mu_r\ |\ r\equiv -i \mod k\}$

We can now define three families of polynomials

\begin{equation}
\tilde G^{(k)}_\lambda(X;q)=
\sum_{T\in\tab_k(\lambda,\, \cdot\,)}q^{\cs(T)}\, x^T
\end{equation}
\begin{equation}
\tH^{(k)}_\mu(X;q)=\sum_{T\in\tab_k(k\mu,\, \cdot\,)}q^{\cs(T)}\, x^T
=\tilde G^{(k)}_{k\mu}(X;q)
\end{equation}
\begin{equation}
H^{(k)}_\mu(X;q)=\sum_{T\in\tab_k(k\mu,\, \cdot\,)}q^{s(T)}\, x^T
=q^{s^*_k(k\mu)}\tH^{(k)}_\mu(X;1/q) \ .
\end{equation}

The parameter $k$ will be called the {\it level} of the corresponding
symmetric functions.

\begin{theorem}\label{Csym}{\rm (symmetry) }
The polynomials $\tG^{(k)}_\lambda$, $\tH^{(k)}_\mu$ and $H^{(k)}_\mu$ are
symmetric.
\end{theorem}

This property follows from Corollary \ref{UV}. Indeed, the commutation
relation $[V_i,V_j]=0$ proves that if $\alpha$ is a rearrangement
of a partition $\mu$,
$$
V_{\alpha_r}\cdots V_{\alpha_1}|0\>=
V_{\mu_r}\cdots V_{\mu_1}|0\> = \psi^k_q(h_\mu)
$$
which shows that for any partition $\lambda$, the sets
$\tab_k(\lambda,\mu)$ and $\tab_k (\lambda,\alpha)$ have the same
spin polynomials.

\begin{remark}{\rm If one defines the linear operator $\varphi_q^k$
as the adjoint of $\psi^k_q$ for the standard scalar product,
the $H$-functions can also be defined by the equation
\begin{equation}
H^{(k)}_\lambda (X; q^{-2}) = \varphi^k_q (s_{k\lambda}) \ .
\end{equation}
}
\end{remark}

There is strong experimental evidence for the following conjectures.

\begin{conjecture}\label{Cpos}{\rm (positivity) }
Their coefficients on the basis of Schur functions are polynomials
with nonnegative integer coefficients.
\end{conjecture}

\begin{conjecture}\label{Cmono}{\rm (monotonicity) }
$H^{(k+1)}_\mu-H^{(k)}_\mu$ is positive on the
Schur basis.
\end{conjecture}

\begin{conjecture}\label{Cplet}{\rm (plethysm) }
When $\mu=\nu^k$, for $\zeta$ a primitive $k$-th root of unity,
$$ H^{(k)}_{\nu^k}(\zeta)=(-1)^{(k-1)|\nu|}\ p_k\circ s_\nu $$
and more generally, when $d|k$ and $\zeta$ is a primitive
$d$-th root of unity,
$$
 H^{(k)}_{\nu^k}(\zeta)=
(-1)^{(d-1)|\nu|k/d}p_d^{k/d}\circ s_{\nu} \ .
$$
Equivalently,
$$
H^{(k)}_{\nu^k}(q) \mod 1-q^k =\sum_{i=0}^{k-1}q^k \ell^{(i)}_k \circ s_\nu
$$
\end{conjecture}

The following statements will be proved in the forthcoming sections.

\begin{theorem}\label{THL}
For $k\ge \ell(\mu)$, $H^{(k)}_\mu$ is equal to the Hall-Littlewood function
$Q'_\mu$.
\end{theorem}

\begin{theorem}\label{IHL}
The difference $Q'_\mu - H^{(2)}_\mu$ is nonnegative on the Schur basis.
\end{theorem}

Taking into account the results of \cite{LLT1,LLT2} and \cite{CL}, this
is sufficient to establish  the conjectures for $k=2$ and $k\ge \ell(\mu)$.

\begin{example}{\rm
{\bf (i) } The $3$-quotient of $\lambda=(3,3,3,2,1)$ is
$((1),(1,1),(1))$ and

\begin{eqnarray*}
\tG_{33321}(q) & = &
  m_{31} + (1+q) m_{22} + (2+2q+q^2) m_{211}\\
&&  + (3 + 5q + 3q^2 + q^3) m_{1111}\\
& = &  s_{31} + q s_{22} + (q+q^2) s_{211} + q^3 s_{1111}\\
\end{eqnarray*}

is a $q$-analogue of the product
$$
s_1 s_{11} s_1 = s_{31} + s_{22} + 2 s_{211} + s_{1111} \ .
$$

{\bf (ii) } The $H$-functions associated to the partition $\lambda=(3,2,1,1)$
are
\begin{eqnarray*}
H_{3211}^{(2)} & = & s_{3211} +q\, s_{322} +q\, s_{331} +q\, s_{4111} \\
	       &   & +(q+q^2)\, s_{421} +q^2\, s_{43} +q^2\, s_{511} + q^3\, s_{52}\\
H_{3211}^{(3)} & = & s_{3211} + q\, s_{322}+(q+q^2)\,s_{331}+q\, s_{4111}\\
               &   & +(q+2q^2)\,s_{421}+(q^2+q^3)\,s_{43} +
(q^2+q^3)\,s_{511}\\
               &   & + 2q^3\, s_{52} + q^4 \, s_{61} \\
H_{3211}^{(4)} & = & s_{3211}+q\,s_{322}+(q+q^2)\,s_{331}+q\,s_{4111}\\
	       &   &
+(q+2q^2+q^3)\,s_{421}+(q^2+q^3+q^4)\,s_{43}+(q^2+q^3+q^4)\,s_{511}\\
               &   & +(2q^3+q^4+q^5)\,s_{52} + (q^4+q^5+q^6)\,s_{61} + q^7\,
s_7 \\
	       & = & Q'_{3211}\\
\end{eqnarray*}
and we see that
$s_{3211} < H_{3211}^{(2)} < H_{3211}^{(3)} < H_{3211}^{(4)} = Q'_{3211}$ .

{\bf (iii)} The plethysms of $s_{21}$ with the cyclic characters $\ell_3^{(i)}$
are given by the reduction modulo $1-q^3$ of
\begin{eqnarray*}
H_{222111}^{(3)}& = &
q^{9}s_{6 3}+  (q+1  )q^{7}s_{6 2 1}+q^{6}s_{6 1 1 1}+
( q+1  )q^{7}s_{5 4}+  (q^{3}+2 q^{2}+2 q+1  )q^{5}s_{5 3 1}\\
&&+  (q^{2}+2 q+1  )q^{5}s_{5 2 2}+  (q^{3}+2 q^{2}+
2 q+1  )q^{4}s_{5 2 1 1}+  (q+1  )q^{4}s_{5 1 1 1 1}\\
&&+ (q^{2}+2 q+1  )q^{5}s_{4 4 1}+  (q^{3}+2 q^{2}+3 q+2)q^{4}s_{4 3 2}
  +  (2 q^{3}+3 q^{2}+3 q+1  )q^{3}s_{ 4 3 1 1}\\
&&+  (q^{3}+3 q^{2}+3 q+2  )q^{3}s_{4 2 2 1}+  (q ^{3}+2 q^{2}+2 q+1
)q^{2}s_{4 2 1 1 1}
+q^{3}s_{4 1 1 1 1 1}+ (q^{3}+1  )q^{3}s_{3 3 3}\\
&&  +  (2 q^{3}+3 q^{2}+2 q+1)q^{2}s_{3 3 2 1}+  (q^{2}+2 q+1  )q^{2}s_{3 3 1 1
1}
+  (q^{2}+2 q+1  )q^{2}s_{3 2 2 2}\\
&&+  (q^{3}+2 q^{2}+2 q+1  )qs_{3 2 2 1 1}+  (q+1  )qs_{3 2 1 1 1 1}
+  (q+ 1  )qs_{2 2 2 2 1}+s_{2 2 2 1 1 1}\\
\end{eqnarray*}
Indeed,
\begin{eqnarray*}
H_{222111}^{(3)}   \mod 1-q^3 &=&
(2  s_{5 2 1 1}+s_{2 2 2 2 1}+s_{3 2 1 1 1 1}+3  s_{4 3 1 1}\\
&& +2  s_{3 2 2 1 1}+s_{5 2 2}+3  s_{4 3 2}+3  s_{3 3 2 1}+s_{3 3 1 1 1}
 +s_{3 2 2 2}+s_{5 1 1 1 1}\\
&&+3  s_{4 2 2 1}+2  s_{5 3 1}+2  s_{4 2 1 1 1} +s_{ 5 4}+s_{6 2 1}+s_{4 4 1}
)q^{2}\\
&& + (2  s_{5 2 1 1}+s_{2 2 2 2 1}+s_{3 2 1 1 1 1}+3  s_{4 3 1 1}
 +2  s_{3 2 2 1 1}+s_{5 2 2}\\
&&+3  s_{4 3 2}+3  s_{3 3 2 1} + s_{3 3 1 1 1}+s_{3 2 2 2}+s_{5 1 1 1 1} \\
&&+3  s_{4 2 2 1}+2  s_{5 3 1} +2  s_{4 2 1 1 1}+s_{5 4}+s_{6 2 1}+s_{4 4 1} )
q \\
&& +2  s_{3 3 1 1 1}+s_{6 3}+s_{6 1 1 1}+2  s_{5 3 1}+2  s_{5 2 2}
 +2  s_ {5 2 1 1}+2  s_{4 4 1}\\
&&+2  s_{4 3 2}+3  s_{4 3 1 1}
 +3  s_{4 2 2 1}+2  s_{4 2 1 1 1}+s_{4 1 1 1 1 1}+2  s_{3 3 3}\\
&& +2  s_{3 3 2 1}+2  s_{3 2 2 2 }+s_{2 2 2 1 1 1}+2  s_{3 2 2 1 1}\\
& =&q^2 \ell_3^{(2)}\circ s_{21}+q\ell_3^{(1)}\circ s_{21}
+\ell_3^{(0)}\circ s_{21} \ .
\end{eqnarray*}

}
\end{example}

\section{The case of dominoes}\label{S7}

For $k=2$, the conjectures can be established by means of the combinatorial
constructions of \cite{CL} and \cite{KLLT}. In this case, conjectures
\ref{Csym}, \ref{Cpos} and \ref{Cplet} follow directly from the results
of \cite{CL}, and the only point remaining to be proved is Theorem \ref{IHL}.

The important special feature of domino tableaux is that there exits a
natural notion of {\it Yamanouchi domino tableau}. These tableaux correspond
to highest weight vectors in tensor products of two irreducible $GL_n$-modules,
in the same way as ordinary Yamanouchi tableaux are the natural labels for
highest weight vectors of irreducible representations.

The {\it column reading} of a domino tableau $T$ is the word obtained by
reading the
successive columns of $T$ from top to bottom and left to right. Horizontal
dominoes, which
belong to two succesive columns $i$ and $i+1$ are read only once, when reading
column $i$.
For example, the column reading of the domino tableau

\begin{center}
\setlength{\unitlength}{0.01in}
\begin{picture}(80,115)(0,-10)
\path(0,60)(20,60)
\path(0,40)(40,40)
\path(20,40)(20,0)
\path(60,40)(60,0)
\path(20,20)(60,20)
\path(0,100)(20,100)(20,60)
	(40,60)(40,40)(80,40)
	(80,0)(0,0)(0,100)
\put(5,5){\makebox(0,0)[lb]{\raisebox{0pt}[0pt][0pt]{\shortstack[l]{{\twlrm
1}}}}}
\put(25,25){\makebox(0,0)[lb]{\raisebox{0pt}[0pt][0pt]{\shortstack[l]{{\twlrm
2}}}}}
\put(5,45){\makebox(0,0)[lb]{\raisebox{0pt}[0pt][0pt]{\shortstack[l]{{\twlrm
3}}}}}
\put(5,85){\makebox(0,0)[lb]{\raisebox{0pt}[0pt][0pt]{\shortstack[l]{{\twlrm
4}}}}}
\put(45,5){\makebox(0,0)[lb]{\raisebox{0pt}[0pt][0pt]{\shortstack[l]{{\twlrm
1}}}}}
\put(65,25){\makebox(0,0)[lb]{\raisebox{0pt}[0pt][0pt]{\shortstack[l]{{\twlrm
2}}}}}
\end{picture}
\end{center}
is $\col(T)=431212$.

A {\it Yamanouchi word} is a word $w=x_1x_2\cdots x_n$ such that each right
factor
$v=x_i\cdots x_n$ of $w$ satisfies $|v|_j\ge |v|_{j+1}$ for each $j$, where
$|v|_j$
denotes the number of occurences of the letter $j$ in $v$.

A {\it Yamanouchi domino tableau} is a domino tableau whose column reading is a
Yamanouchi word. We denote by $\yam_2(\lambda,\mu)$ the set of Yamanouchi
domino
tableaux of shape $\lambda$ and weight $\mu$.

It follows from the results of \cite{CL}, Section 7, that the Schur expansions
of the
$H$-functions of level $2$ are given by
\begin{equation}
H^{(2)}_\lambda =
\sum_\mu \sum_{T\in\yam_2(2\lambda,\mu)}q^{s(T)} s_\mu \ .
\end{equation}
On the other hand,
\begin{equation}
Q'_\lambda =
\sum_\mu \sum_{\t\in\tab(\mu,\lambda)} q^{c(\t)}s_\mu \ .
\end{equation}
To prove Theorem \ref{IHL}, it is thus sufficient to exhibit an injection
$$
\eta :\quad \yam_2(2\lambda,\mu) \longrightarrow  \tab(\mu,\lambda)
$$
satisfying
$$
c(\eta(T)) = s(T) \ .
$$
To achieve this, we shall make use of a bijection described in \cite{BV}, and
extended in \cite{KLLT}, which sends a domino tableau $T\in\tab_2(\alpha,\mu)$
over the alphabet $X=\{1,\ldots,n\}$, to an ordinary tableau
$\t=\phi(T)\in\tab(\alpha,\bar{\mu}\mu)$
over the alphabet $\bar{X}\cup X=\{\bar n<\ldots < \bar 1<1<\ldots <n\}$. The
weight $\bar\mu\mu$
means that $\t$ contains  $\mu_i$ occurences of $i$ and of $\bar i$. The
tableau $\phi(T)$
is invariant under Sch\"utzenberger's involution $\Omega$, and the spin of $T$
can be recovered
from $\t$ by the following procedure \cite{KLLT2}.

Let $\alpha=2\lambda$, $\beta=\alpha'$, $\beta_{\rm
odd}=(\beta_1,\beta_3,\ldots\,)$
and $\beta_{\rm even}=(\beta_2,\beta_4,\ldots\,)$. Then, there exists a
unique factorisation $\t=\tau_1\tau_2$ in the plactic monoid ${\rm
Pl\,}(X\cup\bar X)$,
such that $\tau_1$ is a contretableau of shape $\alpha^1=(\beta_{\rm even})'$
and $\tau_2$ is a tableau of shape $\alpha^2=(\beta_{\rm odd})'$. The
spin of $T=\phi^{-1}(\t)$ is then equal to the number $|\tau_1|_+$ of positive
letters
in $\tau_1$, which is also equal to the number $|\tau_2|_-$ of negative letters
in
$\tau_2$. Moreover, $\tau_2=\Omega(\tau_1)$.

\begin{example}{\rm With the following  tableau $T$ of shape $(4,4,2,2)$, one
finds

\begin{center}
\setlength{\unitlength}{0.01in}
\begingroup\makeatletter\ifx\SetFigFont\undefined
% extract first six characters in \fmtname
\def\x#1#2#3#4#5#6#7\relax{\def\x{#1#2#3#4#5#6}}%
\expandafter\x\fmtname xxxxxx\relax \def\y{splain}%
\ifx\x\y   % LaTeX or SliTeX?
\gdef\SetFigFont#1#2#3{%
  \ifnum #1<17\tiny\else \ifnum #1<20\small\else
  \ifnum #1<24\normalsize\else \ifnum #1<29\large\else
  \ifnum #1<34\Large\else \ifnum #1<41\LARGE\else
     \huge\fi\fi\fi\fi\fi\fi
  \csname #3\endcsname}%
\else
\gdef\SetFigFont#1#2#3{\begingroup
  \count@#1\relax \ifnum 25<\count@\count@25\fi
  \def\x{\endgroup\@setsize\SetFigFont{#2pt}}%
  \expandafter\x
    \csname \romannumeral\the\count@ pt\expandafter\endcsname
    \csname @\romannumeral\the\count@ pt\endcsname
  \csname #3\endcsname}%
\fi
\fi\endgroup
\begin{picture}(480,95)(0,-10)
\path(80,80)(120,80)
\path(120,80)(120,40)(160,40)
	(160,0)(80,0)(80,80)
\path(80,60)(120,60)
\path(80,40)(120,40)
\path(100,40)(100,0)
\path(120,40)(120,0)
\path(120,20)(160,20)
\path(200,40)(300,40)
\path(292.000,38.000)(300.000,40.000)(292.000,42.000)
\path(400,80)(400,0)(480,0)
	(480,40)(440,40)(440,80)(400,80)
\path(400,60)(440,60)
\path(400,40)(440,40)(440,0)
\path(400,20)(480,20)
\path(420,80)(420,0)
\path(460,40)(460,0)
\put(85,5){\makebox(0,0)[lb]{\smash{{{\SetFigFont{12}{14.4}{rm}1}}}}}
\put(105,25){\makebox(0,0)[lb]{\smash{{{\SetFigFont{12}{14.4}{rm}1}}}}}
\put(125,5){\makebox(0,0)[lb]{\smash{{{\SetFigFont{12}{14.4}{rm}1}}}}}
\put(145,25){\makebox(0,0)[lb]{\smash{{{\SetFigFont{12}{14.4}{rm}2}}}}}
\put(85,45){\makebox(0,0)[lb]{\smash{{{\SetFigFont{12}{14.4}{rm}2}}}}}
\put(105,65){\makebox(0,0)[lb]{\smash{{{\SetFigFont{12}{14.4}{rm}3}}}}}
\put(405,5){\makebox(0,0)[lb]{\smash{{{\SetFigFont{12}{14.4}{rm}$\bar 3$}}}}}
\put(425,5){\makebox(0,0)[lb]{\smash{{{\SetFigFont{12}{14.4}{rm}$\bar 2$}}}}}
\put(445,5){\makebox(0,0)[lb]{\smash{{{\SetFigFont{12}{14.4}{rm}$\bar 1$}}}}}
\put(465,5){\makebox(0,0)[lb]{\smash{{{\SetFigFont{12}{14.4}{rm}1}}}}}
\put(465,25){\makebox(0,0)[lb]{\smash{{{\SetFigFont{12}{14.4}{rm}2}}}}}
\put(445,25){\makebox(0,0)[lb]{\smash{{{\SetFigFont{12}{14.4}{rm}1}}}}}
\put(425,25){\makebox(0,0)[lb]{\smash{{{\SetFigFont{12}{14.4}{rm}$\bar 1$}}}}}
\put(405,25){\makebox(0,0)[lb]{\smash{{{\SetFigFont{12}{14.4}{rm}$\bar 2$}}}}}
\put(405,45){\makebox(0,0)[lb]{\smash{{{\SetFigFont{12}{14.4}{rm}$\bar 1$}}}}}
\put(405,65){\makebox(0,0)[lb]{\smash{{{\SetFigFont{12}{14.4}{rm}1}}}}}
\put(425,65){\makebox(0,0)[lb]{\smash{{{\SetFigFont{12}{14.4}{rm}3}}}}}
\put(425,45){\makebox(0,0)[lb]{\smash{{{\SetFigFont{12}{14.4}{rm}2}}}}}
\put(335,40){\makebox(0,0)[lb]{\smash{{{\SetFigFont{12}{14.4}{rm} }}}}}
\put(340,35){\makebox(0,0)[lb]{\smash{{{\SetFigFont{12}{14.4}{rm}${\bf t}$
=}}}}}
\put(0,35){\makebox(0,0)[lb]{\smash{{{\SetFigFont{12}{14.4}{rm}$T$ =}}}}}
\put(245,45){\makebox(0,0)[lb]{\smash{{{\SetFigFont{12}{14.4}{rm}$\phi$}}}}}
\end{picture}
\end{center}
By {\it jeu de taquin}, we find that in the plactic monoid
$$
\t \ \  = \  \
\young{\bar 1 & 1 \cr \bar 2 & \bar 1\cr \blk &\bar 2\cr \blk & \bar 3\cr}
\
\young{3\cr 2\cr 1 & 2 \cr \bar 1 & 1\cr}
\quad =  \ \tau_1\tau_2 \ .
$$
The number of positive letters of $\tau_1$ and the number of
negative letters of $\tau_2$ are both equal to $1$, which is
the spin of $T$.
}
\end{example}

This correspondence still works in the general case ($\alpha$ need not
be of the form $2\lambda$) and the invariant tableau associated to a
domino tableau $T$ admits a similar factorisation $\t=\tau_1\tau_2$,
but in general $\tau_2\not =\Omega(\tau_1)$ and the formula for the
spin is $s(T)={1\over 2}(|\tau_1|_+ +|\tau_2|_-)$.

The map $\eta :\ \yam_2(2\lambda,\mu)\longrightarrow \tab(\mu,\lambda)$
is given by the following algorithm: to compute
$\eta (T)$,
\begin{enumerate}

\item construct the invariant tableau $\t=\phi(T)$

\item apply the {\it jeu de taquin}  algorithm to $\t$ to
obtain the plactic factorization $\t=\tau_1\tau_2$,
and keep only $\tau_2$.

\item Apply the evacuation algorithm to the {\it negative} letters
of $\tau_2$, keeping track of the successive stages.
After all the negative letters have been evacuated, one
is left with a Yamanouchi tableau $\tau$ in positive letters.

\item Complete the tableau $\tau$ to obtain the tableau $\t'=\eta(T)$
using the following rule: suppose that at some stage of the evacuation,
the box of $\tau_2$ which disappeared after the elimination of $\bar i$  was
in row $j$ of $\tau_2$. Then add a box numbered $j$ to row $i$ of $\tau$.
\end{enumerate}

\begin{theorem}\label{etainj}
The above algorithm defines an injection
$$\eta :\ \yam_2(2\lambda,\mu)\longrightarrow \tab(\mu,\lambda)
$$
satifying $c\circ\eta=s$.
\end{theorem}

\begin{corollary}
$H^{(2)}_\lambda \le Q'_\lambda$
\end{corollary}

\begin{example}{\rm Let $T$ be the following Yamanouchi domino tableau, which
is
of shape $2\lambda=(6,4,4,2,2)$, of weight $\mu=(4,3,2)$ and has spin $s(T)=3$
\begin{center}
\setlength{\unitlength}{0.01in}
\begingroup\makeatletter\ifx\SetFigFont\undefined
% extract first six characters in \fmtname
\def\x#1#2#3#4#5#6#7\relax{\def\x{#1#2#3#4#5#6}}%
\expandafter\x\fmtname xxxxxx\relax \def\y{splain}%
\ifx\x\y   % LaTeX or SliTeX?
\gdef\SetFigFont#1#2#3{%
  \ifnum #1<17\tiny\else \ifnum #1<20\small\else
  \ifnum #1<24\normalsize\else \ifnum #1<29\large\else
  \ifnum #1<34\Large\else \ifnum #1<41\LARGE\else
     \huge\fi\fi\fi\fi\fi\fi
  \csname #3\endcsname}%
\else
\gdef\SetFigFont#1#2#3{\begingroup
  \count@#1\relax \ifnum 25<\count@\count@25\fi
  \def\x{\endgroup\@setsize\SetFigFont{#2pt}}%
  \expandafter\x
    \csname \romannumeral\the\count@ pt\expandafter\endcsname
    \csname @\romannumeral\the\count@ pt\endcsname
  \csname #3\endcsname}%
\fi
\fi\endgroup
\begin{picture}(205,115)(0,-10)
\path(85,100)(125,100)(125,60)
	(165,60)(165,20)(205,20)
	(205,0)(85,0)(85,100)
\path(85,60)(125,60)
\path(105,100)(105,60)
\path(125,60)(125,0)
\path(85,40)(125,40)
\path(105,40)(105,0)
\path(125,20)(165,20)(165,0)
\path(145,60)(145,20)
\put(90,5){\makebox(0,0)[lb]{\smash{{{\SetFigFont{12}{14.4}{rm}1}}}}}
\put(115,25){\makebox(0,0)[lb]{\smash{{{\SetFigFont{12}{14.4}{rm}1}}}}}
\put(90,45){\makebox(0,0)[lb]{\smash{{{\SetFigFont{12}{14.4}{rm}2}}}}}
\put(90,85){\makebox(0,0)[lb]{\smash{{{\SetFigFont{12}{14.4}{rm}3}}}}}
\put(110,65){\makebox(0,0)[lb]{\smash{{{\SetFigFont{12}{14.4}{rm}3}}}}}
\put(130,45){\makebox(0,0)[lb]{\smash{{{\SetFigFont{12}{14.4}{rm}2}}}}}
\put(150,25){\makebox(0,0)[lb]{\smash{{{\SetFigFont{12}{14.4}{rm}2}}}}}
\put(130,5){\makebox(0,0)[lb]{\smash{{{\SetFigFont{12}{14.4}{rm}1}}}}}
\put(170,5){\makebox(0,0)[lb]{\smash{{{\SetFigFont{12}{14.4}{rm}1}}}}}
\put(0,40){\makebox(0,0)[lb]{\smash{{{\SetFigFont{12}{14.4}{rm}$T$    =}}}}}
\end{picture}
\end{center}
Then,
$$
\phi(T)\ = \
\young{ 3 & 3\cr 1 & 2 \cr \bar 1 & \bar 1 & 2 & 2 \cr
          \bar 2 & \bar 2 & \bar 1 & 1 \cr
           \bar 3 & \bar 3 & \bar 2 & \bar 1 & 1 & 1\cr}
\qquad \equiv \qquad
\young{ \bar 1 & 1 & 3 \cr \blk & \bar 1 & 2 \cr
        \blk & \bar 2 & \bar 1\cr
        \blk & \blk & \bar 2\cr
        \blk & \blk & \bar 3\cr}
\
\young{ 3 \cr 2 \cr 1 & 2\cr \bar 2 & 1\cr \bar 3 & \bar 1 & 1\cr}
$$
and the succesive stages of the evacuation process are
$$
\matrix{
\young{3\cr 2\cr 1&2\cr \bar 2&1\cr \bar 3&\bar 1&1\cr}
&
\longrightarrow
&
\young{\times\cr 3\cr 2&2\cr 1&1\cr \bar 2&\bar 1&1\cr}
&
\longrightarrow
&
\young{3\cr 2&\times\cr 1&2\cr \bar 1&1&1\cr}
&
\longrightarrow
&
\young{\times\cr 3\cr 2&2\cr 1&1&1\cr}
\cr
& \bar i = \bar 3 & & \bar i =\bar 2 & & \bar i=\bar 1 & \cr
&      j = 5      & &      j = 3     & &      j =    4 & \cr}
$$
so that we find
$$
\eta(T) =
\young{ 3&5\cr 2&2&3\cr 1&1&1&4\cr}
$$
a tableau of shape $\mu=(4,3,2)$, weight $\lambda=(3,2,2,1,1)$ and charge
$c(\t')=3$.

}
\end{example}

\section{The stable case}\label{S8}

As the $Q'$-functions are known to verify all the conjectured properties
of  $H$-functions,
the stable case of the conjectures will be a consequence of
theorem \ref{THL}. This result will be proved by means of
Shimomura's cell decompositions of unipotent varieties.
%%%%%%%%%%%%%%%%%%%%%%%%%%%%%%%%%%%%%%%%%%%%%%%%%%%%%%%%%%%%%%%%

\subsection{Unipotent varieties}

Let $u\in GL(n,\C)$ be a unipotent element,and let
$\F_\nu^u[\C]$ be the variety of $\nu$-flags of $\C^n$
which are fixed by $u$.

It has been shown by N. Shimomura (\cite{Sh1}, see also
 \cite{HSh}) that the   variety $\F_\nu^u[\C]$
admits a cell decomposition, involving only cells of even real dimensions.
More precisely, this cell decomposition is a partition in locally closed
subvarieties, each being algebraically isomorphic to an affine space.
Thus, the odd-dimensional homology groups are zero, and if
$$
\Pi_{\nu\mu}(t^2)=\sum_i t^{2i} \dim H_{2i}(\F_\nu^u,\Z)
$$
is the Poincar\'e polynomial of $\F_\nu^u[\C]$, one has
$|\F_\nu^u[\GF_q] |=\Pi_{\nu\mu}(q)$. But this is also equal to
$\tG_{\nu\mu}(q)$, and as this is true for an infinite set of values
of $q$, one has $\Pi_{\nu\mu}(z)=\tG_{\nu\mu}(z)$ as polynomials.
That is, the coefficient of $\tQ_\mu$ on the monomial function
$m_\nu$ is the Poincar\'e polynomial of $\F_\nu^u$, for a unipotent $u$
of type $\mu$.

Writing
\begin{equation}
\tQ_{\mu} = \sum_{\lambda,\nu}\tK_{\lambda\mu(q)}K_{\lambda\nu}\, m_\nu \ ,
\end{equation}
one sees that
\begin{equation}
\tG_{\nu\mu}(q) =
\sum_{(\t_1,\t_2)\in\tab(\lambda,\mu)\times\tab(\nu,\mu)} q^{\cc(\t_1)} \ .
\end{equation}
Knuth's extension of the Robinson-Schensted correspondence \cite{Kn}
is a bijection between the set
$$
\coprod_\lambda \tab(\lambda,\mu)\times\tab(\lambda,\nu)
$$
of pairs of tableaux with the same shape, and the double coset space
$\S_\mu\backslash \S_n/\S_\nu$ of the symmetric group $\S_n$ modulo
two parabolic subgroups. Double cosets can be encoded by two-line arrays,
integer matrices with prescribed row and column sums, or by {\it tabloids}.

Let $\nu$ and $\mu$ be arbitrary compositions of the same integer $n$.
A $\mu$-tabloid of shape $\nu$ is a filling of the diagram of boxes
with row lengths $\nu_1,\nu_2,\ldots,\nu_r$,  the lowest row being
numbered $1$ (French convention for tableaux), such that the number $i$
occurs $\mu_i$ times, and such that each row is nondecreasing. For example,
$$
\young{ 3 \cr 1 & 1 & 1 \cr 1 & 1& 3 \cr 2 & 3 \cr}
$$
is a $(5,1,3)$-tabloid of shape $(2,3,3,1)$.

We denote by $L(\nu,\mu)$ the set of tabloids of shape $\nu$ and
weight $\mu$.  A tabloid will be identified with the word obtained by reading
it from left to right and top to bottom.
Then,
\begin{equation}
\tG_{\nu\mu}(q) =\sum_{T\in L(\nu,\mu)}q^{\cc(T)} \ .
\end{equation}
\begin{example}{\rm To compute $\tG_{42,321}(q)$ one lists the elements
of $L((4,2),(3,2,1))$, which are
$$
\young{2&3\cr 1&1&1&2\cr}\qquad
\young{2&2\cr 1&1&1&3\cr}\qquad
\young{1&3\cr 1&1&2&2\cr}\qquad
\young{1&2\cr 1&1&2&3\cr}\qquad
\young{1&1\cr 1&2&2&3\cr}
$$
Reading them as prescribed,
we obtain the words
$$
231112\qquad 221113\qquad 131122\qquad 121123\qquad 111223
$$
whose respective charges are $2,1,3,2,4$. The cocharge polynomial is
thus $\tG_{42,321}(q) = 1+q+2q^2+q^3$.
}
\end{example}

In Shimomura's decomposition of the fixed point variety $\F_\mu^u$ of a
unipotent
of type $\nu$, the cells are indexed by tabloids of shape $\nu$ and weight
$\mu$. The dimension $d(T)$ of the cell $c_T$ indexed by $T\in L(\nu,\mu)$ is
computed by an algorithm described below, and gives another combinatorial
interpretation
of the polynomial $\tG_{\mu\nu}(q)$, exchanging the r\^oles of shape and
weight:

\begin{equation}
\tG_{\mu\nu}=\sum_{T\in L(\mu,\nu)}q^{\cc(T)}
=\sum_{T\in L(\nu,\mu)}q^{d(T)} \ .
\end{equation}
The dimensions $d(T)$ are given by the following algorithm.

\begin{enumerate}

\item If $T\in L(\nu,(n))$ then $d(T)=0$;

\item If $\mu=(\mu_1,\mu_2)$ has exactly two parts, and $T\in L(\nu,\mu)$,
then $d(T)$ is computed as follows. A box $\alpha$ of $T$ is said to
be {\it special} if it contains the rightmost $1$ of its row. For a box
$\beta$ of $T$, put $d(\beta)$=0 if $\beta$ does not contain a $2$,
and if $\beta$ contains a $2$, set $d(\beta)$ equal to the number of
nonspecial $1$'s lying in the column of $\beta$, plus the number of
special $1$'s lying in the same column, but in a lower position. Then
$$d(T)=\sum_\beta d(\beta)\ .$$

\item Let $\mu=(\mu_1,\ldots,\mu_k)$ and $\mu^*=(\mu_1,\ldots,\mu_{k-1})$.
For $T\in L(\nu,\mu)$, let $T_1$ be the tabloid obtained by changing
the entries $k$ into $2$ and all the other ones into $1$. Let $T_2$
be the tabloid obtained by erasing all the entries $k$, {\it and rearranging
the rows in the appropriate order}. Then,
\begin{equation}
d(T)=d(T_1)+d(T_2) \ .
\end{equation}
\end{enumerate}

\begin{example}{\rm
With $T=\young{1&4\cr1&2&3\cr1&1&2\cr} \in L(332,4211)$, one has
$$
T_1:=\young{1&2\cr1&1&{\bf 1}\cr1&1&{\bf 1}\cr}
\qquad
T_2=\young{1\cr 1&2&3 \cr 1&1&2\cr}
\qquad
T_{21}=\young{1\cr 1&{\bf 1}&2\cr 1&1&{\bf 1}\cr}
\qquad
T_{22}=\young{1\cr {\bf 1}&2\cr 1&{\bf 1}&2\cr}
$$
where the special entries are printed in boldtype.
Thus, $d(T)=t(T_1)+d(T_2)=2+d(T_{21})+d(T_{22})=4$.
}
\end{example}

We shall need a variant of this construction, in which the shape
$\nu$ is allowed  to be an arbitrary composition, and where in step 3,
the rearranging of the rows is supressed. Such a variant has already
been used by Terada \cite{Te} in the case of complete flags.

That is, we associate to a tabloid $T\in L(\nu,\mu)$ an integer $e(T)$,
defined by

\begin{enumerate}

\item For $T\in L(\nu,(n))$, $e(T)=d(T)=0$;

\item For $T\in L(\nu, (\mu_1,\mu_2))$, $e(T)=d(T)$;

\item Otherwise $e(T)=e(T_1)+e(T_2)$ where $T_1$ is defined as above,
but this time $T_2$ is obtained from $T$ by erasing the entries $k$,
without reordering.
\end{enumerate}

\begin{lemma}
Let $\lambda=(\lambda_1,\ldots,\lambda_r)$ be a partition, and let
$\nu=\lambda\cdot\sigma=(\lambda_{\sigma(1)},\ldots,\lambda_{\sigma(r)})$,
$\sigma\in \S_r$. Then,
the distribution of $e$ on $L(\nu,\mu)$ is the same as the distribution of
$d$ on $L(\lambda,\mu)$.
That is,
$$D_{\lambda\mu}(q)=\sum_{T\in L(\lambda,\mu)}q^{d(T)}
=E_{\nu\mu}(q)=\sum_{T\in L(\nu,\mu)}q^{e(T)} \ .
$$
In particular, $D_{\lambda\mu}(q)=E_{\lambda\mu}(q)$.
\end{lemma}

\Proof This could be proved by repeating word for word the geometric argument
of \cite{Sh1}. We give here a short combinatorial argument. As the two
statistics
coincide on tabloids whose shape is a partition and whose weight has at most
two parts, the only thing to prove, thanks to the recurrence formula, is that
$e$ has the same distribution
on $L(\beta,(\mu_1,\mu_2))$ as on $L(\alpha,(\mu_1,\mu_2))$ when $\beta$ is
a permutation of $\alpha$. The symmetric group being generated by the
elementary
transpositions $\sigma_i=(i,i+1)$, one may assume that $\beta=\alpha\sigma_i$.
We define the image $T\sigma_i$ of a tabloid $T\in L(\alpha,(\mu_1,\mu_2))$ by
distinguishing among the following configurations for rows $i$ and $i+1$:

\begin{enumerate}

\item
$$
\left.\matrix{ x_1 &\ldots & x_k&2&2^r\cr
		1  &\ldots & 1 &{\bf 1} & 2^s\cr}\right.
	\qquad
{\sigma_i \atop \makebox[1cm]{\rightarrowfill} }
\qquad
\left.\matrix{ x_1 &\ldots & x_k&2&2^s\cr
               1  &\ldots & 1 &{\bf 1} & 2^r\cr}\right.
$$
\item
$$
\left.\matrix{1  &\ldots & 1 &{\bf 1} & 2^r\cr
              x_1 &\ldots & x_k&2&2^s\cr}\right.
	      \qquad
	      {\sigma_i\atop \makebox[1cm]{\rightarrowfill}}
	      \qquad
\left.\matrix{1  &\ldots & 1 &{\bf 1} & 2^s\cr
	      x_1 &\ldots & x_k&2&2^r\cr}\right.
$$
\item In all other cases, the two rows are exchanged:
$$
\left.\matrix { x_1&\ldots & x_r\cr y_1&\ldots &&y_s\cr}\right.
\qquad
{\sigma_i\atop \makebox[1cm]{\rightarrowfill}}
\qquad
\left.\matrix {y_1&\ldots &&y_s\cr  x_1&\ldots & x_r\cr}\right.
$$
\end{enumerate}
{}From this definition, it is clear that $e(T\sigma_i)=e(T)$. Moreover, it
is not difficult to check that this defines an $e$-preserving action
of the symmetric group $\S_m$ on the set of $\mu$-tabloids with $m$
rows, such that $L(\alpha,\mu)\sigma=L(\alpha\sigma,\mu)$ (the only point
needing a verification is the braid relation
$\sigma_i\sigma_{i+1}\sigma_i=\sigma_{i+1}\sigma_i\sigma_{i+1}$).

Thus, for a partition $\lambda$ and a two-part weight $\mu=(\mu_1,\mu_2)$,
$d$ and $e$ coincide on $L(\lambda,\mu)$, and for $\sigma\in\S_m$,
$E_{\lambda\sigma,\mu}(q)=D_{\lambda\mu}(q)$. Now, by induction, for
$\mu=(\mu_1,\ldots,\mu_k)$,
$$
D_{\lambda\mu}(q)=\sum_{T\in L(\lambda,\mu)}q^{d(T_1)}q^{d(T_2)}
$$
$$
=\sum_{\bar\lambda={\rm shape\,}(T_1)}q^{d(T_1)}D_{\bar\lambda,\mu^*}(q)
=\sum_{\bar\lambda={\rm shape\,}(T_1)}e^{e(T_1)}E_{\bar\lambda,\mu^*}(q)
=E_{\lambda\mu}(q) \ .
$$
\cqfd

\begin{example}{\rm
Take $\lambda=(3,2,1)$, $\mu=(4,2)$ and $\nu=\lambda\sigma_1\sigma_2=(3,1,2)$.
The $\mu$-tabloids of shape $\lambda$ are
$$
\matrix{
T
&
\young{2\cr 1&2\cr 1&1&1\cr}
&
\young{2\cr 1&1\cr 1&1&2\cr}
&
\young{1\cr 1&1\cr 1&2&2\cr}
&
\young{1\cr 2&2\cr 1&1&1\cr}
&
\young{1\cr 1&2\cr 1&1&2\cr}
\cr
d(T) & 3 & 2 & 0 & 2 & 1 \cr}
$$
The $\nu$-tabloids of shape $\lambda$ are
$$
\matrix{
T
&
\young{1&1&1\cr 1\cr 2&2\cr}
&
\young{1&1&1\cr 2\cr 1&2\cr}
&
\young{1&1&2\cr 1\cr 1&2\cr}
&
\young{1&1&2\cr 2\cr 1&1\cr}
&
\young{1&2&2\cr 1\cr 1&1\cr}
\cr
e(T) & 2 & 3 & 0 & 2 & 1\cr}
$$
Thus, $D_{\lambda\mu}(q)=E_{\nu\mu}(q)=1+q+2q^2+q^3=\tG_{\mu\lambda}(q)$.
The tabloids contributing a term $q^2$ are apparied in the following way:
$$
\young{2\cr 1&1\cr 1&1&2\cr}\quad\longrightarrow\quad
\young{1&1&2\cr 2\cr 1&1\cr}\qquad
\young{1\cr 2&2\cr 1&1&1\cr}\quad\longrightarrow\quad
\young{1&1&1\cr 1\cr 2&2\cr}
$$
}
\end{example}

\begin{remark}{\rm
 The only property that we shall need in the sequel is the
equality $D_{\lambda\mu}(q)=E_{\lambda\mu}(q)$. However, it is
possible to be more explicit by constructing a bijection exchanging
$d$ and $e$. The above action of $\S_m$ can be extended to tabloids
with arbitrary weight, still preserving $e$. Suppose for example
that we want to apply $\sigma_i$ to a tabloid $T$ whose restriction
to rows $i,i+1$ is
$$
\young{1&1&2&3&7&7&9\cr
       1&1&1&2&6&6&6&8&8&9\cr}
$$
One first determines the positions of the greatest entries, which
are the $9$'s, in $T\sigma_i$. Starting with an empty diagram
of the permuted shape $(10,7)$, one constructs $T_1$ as above
by converting all the entries $9$ of $T$ into $2$ and the remaining
ones into $1$. Then we apply $\sigma_i$ to $T_1$, and the positions
of the $2$ in $T_1\sigma_i$ give the positions of the $9$ in  $T\sigma_i$.
Then, the entries $9$ are removed from $T$ ad the procedure is iterated
until one reaches a tabloid whose rows $i$ and $i+1$ are of equal lenghts.
This tabloid is then copied (without permutation) in the remaining part
of the result. On the example, this gives

\small
\begin{center}
\setlength{\unitlength}{0.008in}
\begin{picture}(760,455)(0,-10)
\path(580,440)(580,400)
\path(600,440)(600,400)
\path(620,440)(620,400)
\path(640,440)(640,400)
\path(660,440)(660,400)
\path(680,440)(680,400)
\path(700,440)(700,420)
\path(720,440)(720,420)
\path(740,440)(740,420)
\path(560,420)(700,420)
\path(560,440)(560,400)(700,400)
	(700,420)(760,420)(760,440)(560,440)
\path(560,340)(560,300)(700,300)
	(700,320)(760,320)(760,340)(560,340)
\path(580,340)(580,300)
\path(600,340)(600,300)
\path(620,340)(620,300)
\path(640,340)(640,300)
\path(660,340)(660,300)
\path(680,340)(680,300)
\path(700,340)(700,320)
\path(720,340)(720,320)
\path(740,340)(740,320)
\path(560,320)(700,320)
\path(560,340)(560,300)(700,300)
	(700,320)(760,320)(760,340)(560,340)
\path(560,240)(560,200)(700,200)
	(700,220)(760,220)(760,240)(560,240)
\path(580,240)(580,200)
\path(600,240)(600,200)
\path(620,240)(620,200)
\path(640,240)(640,200)
\path(660,240)(660,200)
\path(680,240)(680,200)
\path(700,240)(700,220)
\path(720,240)(720,220)
\path(740,240)(740,220)
\path(560,220)(700,220)
\path(560,240)(560,200)(700,200)
	(700,220)(760,220)(760,240)(560,240)
\path(560,140)(560,100)(700,100)
	(700,120)(760,120)(760,140)(560,140)
\path(580,140)(580,100)
\path(600,140)(600,100)
\path(620,140)(620,100)
\path(640,140)(640,100)
\path(660,140)(660,100)
\path(680,140)(680,100)
\path(700,140)(700,120)
\path(720,140)(720,120)
\path(740,140)(740,120)
\path(560,120)(700,120)
\path(560,140)(560,100)(700,100)
	(700,120)(760,120)(760,140)(560,140)
\path(0,440)(140,440)(140,420)
	(200,420)(200,400)(0,400)(0,440)
\path(0,420)(140,420)
\path(20,440)(20,400)
\path(40,440)(40,400)
\path(60,440)(60,400)
\path(80,440)(80,400)
\path(100,440)(100,400)
\path(120,440)(120,400)
\path(140,420)(140,400)
\path(160,420)(160,400)
\path(180,420)(180,400)
\path(220,420)(280,420)
\path(272.000,418.000)(280.000,420.000)(272.000,422.000)
\path(300,440)(300,400)(440,400)
	(440,420)(500,420)(500,440)(300,440)
\path(320,440)(320,400)
\path(340,440)(340,400)
\path(360,440)(360,400)
\path(380,440)(380,400)
\path(400,440)(400,400)
\path(420,440)(420,400)
\path(440,440)(440,420)
\path(460,440)(460,420)
\path(480,440)(480,420)
\path(300,420)(440,420)
\path(0,340)(0,300)(180,300)
	(180,320)(120,320)(120,340)(0,340)
\path(220,320)(280,320)
\path(272.000,318.000)(280.000,320.000)(272.000,322.000)
\path(300,340)(480,340)(480,320)
	(420,320)(420,300)(300,300)(300,340)
\path(0,240)(0,200)(140,200)
	(140,220)(120,220)(120,240)(0,240)
\path(220,220)(280,220)
\path(272.000,218.000)(280.000,220.000)(272.000,222.000)
\path(300,240)(300,200)(420,200)
	(420,220)(440,220)(440,240)(300,240)
\path(0,140)(0,100)(140,100)
	(140,120)(80,120)(80,140)(0,140)
\path(220,120)(280,120)
\path(272.000,118.000)(280.000,120.000)(272.000,122.000)
\path(300,140)(380,140)(440,140)
	(440,120)(380,120)(380,100)
	(300,100)(300,140)
\path(0,40)(0,0)(80,0)
	(80,40)(0,40)
\path(220,20)(280,20)
\path(272.000,18.000)(280.000,20.000)(272.000,22.000)
\path(300,40)(300,0)(380,0)
	(380,40)(300,40)
\path(560,40)(560,0)(700,0)
	(700,20)(760,20)(760,40)(580,40)
\path(580,40)(560,40)
\path(580,40)(580,0)
\path(600,40)(600,0)
\path(620,40)(620,0)
\path(640,40)(640,0)
\path(660,40)(660,0)
\path(680,40)(680,0)
\path(700,40)(700,20)
\path(720,40)(720,20)
\path(740,40)(740,20)
\path(560,20)(700,20)
\path(0,320)(120,320)
\path(300,320)(420,320)
\path(0,220)(120,220)
\path(300,220)(420,220)
\path(0,120)(80,120)
\path(300,120)(380,120)
\path(0,20)(80,20)
\path(300,20)(380,20)
\path(320,40)(320,0)
\path(340,40)(340,0)
\path(360,40)(360,0)
\path(320,140)(320,100)
\path(340,140)(340,100)
\path(360,120)(360,100)
\path(360,140)(360,120)
\path(380,140)(380,120)
\path(400,140)(400,120)
\path(420,140)(420,120)
\path(320,240)(320,200)
\path(340,240)(340,200)
\path(360,240)(360,200)
\path(380,240)(380,200)
\path(400,240)(400,200)
\path(420,240)(420,220)
\path(320,340)(320,300)
\path(340,340)(340,300)
\path(360,340)(360,300)
\path(380,340)(380,300)
\path(400,340)(400,300)
\path(420,340)(420,320)
\path(440,340)(440,320)
\path(460,340)(460,320)
\path(20,340)(20,300)
\path(40,340)(40,300)
\path(60,340)(60,300)
\path(80,340)(80,300)
\path(100,340)(100,300)
\path(120,320)(120,300)
\path(140,320)(140,300)
\path(160,320)(160,300)
\path(20,240)(20,200)
\path(40,240)(40,200)
\path(60,240)(60,200)
\path(80,240)(80,200)
\path(100,240)(100,200)
\path(120,220)(120,200)
\path(20,140)(20,100)
\path(40,140)(40,100)
\path(60,140)(60,100)
\path(80,120)(80,100)
\drawline(100,120)(100,120)
\path(100,120)(100,100)
\path(120,120)(120,100)
\path(20,40)(20,0)
\path(40,40)(40,0)
\path(60,40)(60,0)
\put(5,425){\makebox(0,0)[lb]{\raisebox{0pt}[0pt][0pt]{\shortstack[l]{{\twlrm
1}}}}}
\put(25,425){\makebox(0,0)[lb]{\raisebox{0pt}[0pt][0pt]{\shortstack[l]{{\twlrm
1}}}}}
\put(45,425){\makebox(0,0)[lb]{\raisebox{0pt}[0pt][0pt]{\shortstack[l]{{\twlrm
1}}}}}
\put(65,425){\makebox(0,0)[lb]{\raisebox{0pt}[0pt][0pt]{\shortstack[l]{{\twlrm
1}}}}}
\put(85,425){\makebox(0,0)[lb]{\raisebox{0pt}[0pt][0pt]{\shortstack[l]{{\twlrm
1}}}}}
\put(105,425){\makebox(0,0)[lb]{\raisebox{0pt}[0pt][0pt]{\shortstack[l]{{\twlrm
1}}}}}
\put(125,425){\makebox(0,0)[lb]{\raisebox{0pt}[0pt][0pt]{\shortstack[l]{{\twlrm
2}}}}}
\path(560,440)(560,400)(700,400)
	(700,420)(760,420)(760,440)(560,440)
\put(5,405){\makebox(0,0)[lb]{\raisebox{0pt}[0pt][0pt]{\shortstack[l]{{\twlrm
1}}}}}
\put(685,5){\makebox(0,0)[lb]{\raisebox{0pt}[0pt][0pt]{\shortstack[l]{{\twlrm
9}}}}}
\put(25,405){\makebox(0,0)[lb]{\raisebox{0pt}[0pt][0pt]{\shortstack[l]{{\twlrm
1}}}}}
\put(45,405){\makebox(0,0)[lb]{\raisebox{0pt}[0pt][0pt]{\shortstack[l]{{\twlrm
1}}}}}
\put(65,405){\makebox(0,0)[lb]{\raisebox{0pt}[0pt][0pt]{\shortstack[l]{{\twlrm
1}}}}}
\put(85,405){\makebox(0,0)[lb]{\raisebox{0pt}[0pt][0pt]{\shortstack[l]{{\twlrm
1}}}}}
\put(105,405){\makebox(0,0)[lb]{\raisebox{0pt}[0pt][0pt]{\shortstack[l]{{\twlrm
1}}}}}
\put(125,405){\makebox(0,0)[lb]{\raisebox{0pt}[0pt][0pt]{\shortstack[l]{{\twlrm
1}}}}}
\put(145,405){\makebox(0,0)[lb]{\raisebox{0pt}[0pt][0pt]{\shortstack[l]{{\twlrm
1}}}}}
\put(165,405){\makebox(0,0)[lb]{\raisebox{0pt}[0pt][0pt]{\shortstack[l]{{\twlrm
1}}}}}
\put(185,405){\makebox(0,0)[lb]{\raisebox{0pt}[0pt][0pt]{\shortstack[l]{{\twlrm
2}}}}}
\put(305,405){\makebox(0,0)[lb]{\raisebox{0pt}[0pt][0pt]{\shortstack[l]{{\twlrm
1}}}}}
\put(325,405){\makebox(0,0)[lb]{\raisebox{0pt}[0pt][0pt]{\shortstack[l]{{\twlrm
1}}}}}
\put(345,405){\makebox(0,0)[lb]{\raisebox{0pt}[0pt][0pt]{\shortstack[l]{{\twlrm
1}}}}}
\put(365,405){\makebox(0,0)[lb]{\raisebox{0pt}[0pt][0pt]{\shortstack[l]{{\twlrm
1}}}}}
\put(385,405){\makebox(0,0)[lb]{\raisebox{0pt}[0pt][0pt]{\shortstack[l]{{\twlrm
1}}}}}
\put(405,405){\makebox(0,0)[lb]{\raisebox{0pt}[0pt][0pt]{\shortstack[l]{{\twlrm
1}}}}}
\put(425,405){\makebox(0,0)[lb]{\raisebox{0pt}[0pt][0pt]{\shortstack[l]{{\twlrm
2}}}}}
\put(305,425){\makebox(0,0)[lb]{\raisebox{0pt}[0pt][0pt]{\shortstack[l]{{\twlrm
1}}}}}
\put(325,425){\makebox(0,0)[lb]{\raisebox{0pt}[0pt][0pt]{\shortstack[l]{{\twlrm
1}}}}}
\put(345,425){\makebox(0,0)[lb]{\raisebox{0pt}[0pt][0pt]{\shortstack[l]{{\twlrm
1}}}}}
\put(365,425){\makebox(0,0)[lb]{\raisebox{0pt}[0pt][0pt]{\shortstack[l]{{\twlrm
1}}}}}
\put(385,425){\makebox(0,0)[lb]{\raisebox{0pt}[0pt][0pt]{\shortstack[l]{{\twlrm
1}}}}}
\put(405,425){\makebox(0,0)[lb]{\raisebox{0pt}[0pt][0pt]{\shortstack[l]{{\twlrm
1}}}}}
\put(425,425){\makebox(0,0)[lb]{\raisebox{0pt}[0pt][0pt]{\shortstack[l]{{\twlrm
1}}}}}
\put(445,425){\makebox(0,0)[lb]{\raisebox{0pt}[0pt][0pt]{\shortstack[l]{{\twlrm
1}}}}}
\put(465,425){\makebox(0,0)[lb]{\raisebox{0pt}[0pt][0pt]{\shortstack[l]{{\twlrm
1}}}}}
\put(485,425){\makebox(0,0)[lb]{\raisebox{0pt}[0pt][0pt]{\shortstack[l]{{\twlrm
2}}}}}
\put(745,425){\makebox(0,0)[lb]{\raisebox{0pt}[0pt][0pt]{\shortstack[l]{{\twlrm
9}}}}}
\put(685,405){\makebox(0,0)[lb]{\raisebox{0pt}[0pt][0pt]{\shortstack[l]{{\twlrm
9}}}}}
\put(145,305){\makebox(0,0)[lb]{\raisebox{0pt}[0pt][0pt]{\shortstack[l]{{\twlrm
2}}}}}
\put(165,305){\makebox(0,0)[lb]{\raisebox{0pt}[0pt][0pt]{\shortstack[l]{{\twlrm
2}}}}}
\put(5,325){\makebox(0,0)[lb]{\raisebox{0pt}[0pt][0pt]{\shortstack[l]{{\twlrm
1}}}}}
\put(25,325){\makebox(0,0)[lb]{\raisebox{0pt}[0pt][0pt]{\shortstack[l]{{\twlrm
1}}}}}
\put(45,325){\makebox(0,0)[lb]{\raisebox{0pt}[0pt][0pt]{\shortstack[l]{{\twlrm
1}}}}}
\put(65,325){\makebox(0,0)[lb]{\raisebox{0pt}[0pt][0pt]{\shortstack[l]{{\twlrm
1}}}}}
\put(85,325){\makebox(0,0)[lb]{\raisebox{0pt}[0pt][0pt]{\shortstack[l]{{\twlrm
1}}}}}
\put(105,325){\makebox(0,0)[lb]{\raisebox{0pt}[0pt][0pt]{\shortstack[l]{{\twlrm
1}}}}}
\put(5,305){\makebox(0,0)[lb]{\raisebox{0pt}[0pt][0pt]{\shortstack[l]{{\twlrm
1}}}}}
\put(25,305){\makebox(0,0)[lb]{\raisebox{0pt}[0pt][0pt]{\shortstack[l]{{\twlrm
1}}}}}
\put(45,305){\makebox(0,0)[lb]{\raisebox{0pt}[0pt][0pt]{\shortstack[l]{{\twlrm
1}}}}}
\put(65,305){\makebox(0,0)[lb]{\raisebox{0pt}[0pt][0pt]{\shortstack[l]{{\twlrm
1}}}}}
\put(85,305){\makebox(0,0)[lb]{\raisebox{0pt}[0pt][0pt]{\shortstack[l]{{\twlrm
1}}}}}
\put(105,305){\makebox(0,0)[lb]{\raisebox{0pt}[0pt][0pt]{\shortstack[l]{{\twlrm
1}}}}}
\put(125,305){\makebox(0,0)[lb]{\raisebox{0pt}[0pt][0pt]{\shortstack[l]{{\twlrm
1}}}}}
\put(305,325){\makebox(0,0)[lb]{\raisebox{0pt}[0pt][0pt]{\shortstack[l]{{\twlrm
1}}}}}
\put(325,325){\makebox(0,0)[lb]{\raisebox{0pt}[0pt][0pt]{\shortstack[l]{{\twlrm
1}}}}}
\put(345,325){\makebox(0,0)[lb]{\raisebox{0pt}[0pt][0pt]{\shortstack[l]{{\twlrm
1}}}}}
\put(365,325){\makebox(0,0)[lb]{\raisebox{0pt}[0pt][0pt]{\shortstack[l]{{\twlrm
1}}}}}
\put(385,325){\makebox(0,0)[lb]{\raisebox{0pt}[0pt][0pt]{\shortstack[l]{{\twlrm
1}}}}}
\put(405,325){\makebox(0,0)[lb]{\raisebox{0pt}[0pt][0pt]{\shortstack[l]{{\twlrm
1}}}}}
\put(425,325){\makebox(0,0)[lb]{\raisebox{0pt}[0pt][0pt]{\shortstack[l]{{\twlrm
1}}}}}
\put(445,325){\makebox(0,0)[lb]{\raisebox{0pt}[0pt][0pt]{\shortstack[l]{{\twlrm
2}}}}}
\put(465,325){\makebox(0,0)[lb]{\raisebox{0pt}[0pt][0pt]{\shortstack[l]{{\twlrm
2}}}}}
\put(305,305){\makebox(0,0)[lb]{\raisebox{0pt}[0pt][0pt]{\shortstack[l]{{\twlrm
1}}}}}
\put(325,305){\makebox(0,0)[lb]{\raisebox{0pt}[0pt][0pt]{\shortstack[l]{{\twlrm
1}}}}}
\put(345,305){\makebox(0,0)[lb]{\raisebox{0pt}[0pt][0pt]{\shortstack[l]{{\twlrm
1}}}}}
\put(365,305){\makebox(0,0)[lb]{\raisebox{0pt}[0pt][0pt]{\shortstack[l]{{\twlrm
1}}}}}
\put(385,305){\makebox(0,0)[lb]{\raisebox{0pt}[0pt][0pt]{\shortstack[l]{{\twlrm
1}}}}}
\put(405,305){\makebox(0,0)[lb]{\raisebox{0pt}[0pt][0pt]{\shortstack[l]{{\twlrm
1}}}}}
\put(745,325){\makebox(0,0)[lb]{\raisebox{0pt}[0pt][0pt]{\shortstack[l]{{\twlrm
9}}}}}
\put(685,305){\makebox(0,0)[lb]{\raisebox{0pt}[0pt][0pt]{\shortstack[l]{{\twlrm
9}}}}}
\put(705,325){\makebox(0,0)[lb]{\raisebox{0pt}[0pt][0pt]{\shortstack[l]{{\twlrm
8}}}}}
\put(725,325){\makebox(0,0)[lb]{\raisebox{0pt}[0pt][0pt]{\shortstack[l]{{\twlrm
8}}}}}
\put(5,225){\makebox(0,0)[lb]{\raisebox{0pt}[0pt][0pt]{\shortstack[l]{{\twlrm
1}}}}}
\put(25,225){\makebox(0,0)[lb]{\raisebox{0pt}[0pt][0pt]{\shortstack[l]{{\twlrm
1}}}}}
\put(45,225){\makebox(0,0)[lb]{\raisebox{0pt}[0pt][0pt]{\shortstack[l]{{\twlrm
1}}}}}
\put(65,225){\makebox(0,0)[lb]{\raisebox{0pt}[0pt][0pt]{\shortstack[l]{{\twlrm
1}}}}}
\put(85,225){\makebox(0,0)[lb]{\raisebox{0pt}[0pt][0pt]{\shortstack[l]{{\twlrm
2}}}}}
\put(105,225){\makebox(0,0)[lb]{\raisebox{0pt}[0pt][0pt]{\shortstack[l]{{\twlrm
2}}}}}
\put(5,205){\makebox(0,0)[lb]{\raisebox{0pt}[0pt][0pt]{\shortstack[l]{{\twlrm
1}}}}}
\put(25,205){\makebox(0,0)[lb]{\raisebox{0pt}[0pt][0pt]{\shortstack[l]{{\twlrm
1}}}}}
\put(45,205){\makebox(0,0)[lb]{\raisebox{0pt}[0pt][0pt]{\shortstack[l]{{\twlrm
1}}}}}
\put(65,205){\makebox(0,0)[lb]{\raisebox{0pt}[0pt][0pt]{\shortstack[l]{{\twlrm
1}}}}}
\put(85,205){\makebox(0,0)[lb]{\raisebox{0pt}[0pt][0pt]{\shortstack[l]{{\twlrm
1}}}}}
\put(105,205){\makebox(0,0)[lb]{\raisebox{0pt}[0pt][0pt]{\shortstack[l]{{\twlrm
1}}}}}
\put(125,205){\makebox(0,0)[lb]{\raisebox{0pt}[0pt][0pt]{\shortstack[l]{{\twlrm
1}}}}}
\put(305,225){\makebox(0,0)[lb]{\raisebox{0pt}[0pt][0pt]{\shortstack[l]{{\twlrm
1}}}}}
\put(325,225){\makebox(0,0)[lb]{\raisebox{0pt}[0pt][0pt]{\shortstack[l]{{\twlrm
1}}}}}
\put(345,225){\makebox(0,0)[lb]{\raisebox{0pt}[0pt][0pt]{\shortstack[l]{{\twlrm
1}}}}}
\put(365,225){\makebox(0,0)[lb]{\raisebox{0pt}[0pt][0pt]{\shortstack[l]{{\twlrm
1}}}}}
\put(385,225){\makebox(0,0)[lb]{\raisebox{0pt}[0pt][0pt]{\shortstack[l]{{\twlrm
1}}}}}
\put(405,225){\makebox(0,0)[lb]{\raisebox{0pt}[0pt][0pt]{\shortstack[l]{{\twlrm
1}}}}}
\put(425,225){\makebox(0,0)[lb]{\raisebox{0pt}[0pt][0pt]{\shortstack[l]{{\twlrm
1}}}}}
\put(305,205){\makebox(0,0)[lb]{\raisebox{0pt}[0pt][0pt]{\shortstack[l]{{\twlrm
1}}}}}
\put(325,205){\makebox(0,0)[lb]{\raisebox{0pt}[0pt][0pt]{\shortstack[l]{{\twlrm
1}}}}}
\put(345,205){\makebox(0,0)[lb]{\raisebox{0pt}[0pt][0pt]{\shortstack[l]{{\twlrm
1}}}}}
\put(365,205){\makebox(0,0)[lb]{\raisebox{0pt}[0pt][0pt]{\shortstack[l]{{\twlrm
1}}}}}
\put(385,205){\makebox(0,0)[lb]{\raisebox{0pt}[0pt][0pt]{\shortstack[l]{{\twlrm
2}}}}}
\put(405,205){\makebox(0,0)[lb]{\raisebox{0pt}[0pt][0pt]{\shortstack[l]{{\twlrm
2}}}}}
\put(645,205){\makebox(0,0)[lb]{\raisebox{0pt}[0pt][0pt]{\shortstack[l]{{\twlrm
7}}}}}
\put(665,205){\makebox(0,0)[lb]{\raisebox{0pt}[0pt][0pt]{\shortstack[l]{{\twlrm
7}}}}}
\put(685,205){\makebox(0,0)[lb]{\raisebox{0pt}[0pt][0pt]{\shortstack[l]{{\twlrm
9}}}}}
\put(705,225){\makebox(0,0)[lb]{\raisebox{0pt}[0pt][0pt]{\shortstack[l]{{\twlrm
8}}}}}
\put(725,225){\makebox(0,0)[lb]{\raisebox{0pt}[0pt][0pt]{\shortstack[l]{{\twlrm
8}}}}}
\put(745,225){\makebox(0,0)[lb]{\raisebox{0pt}[0pt][0pt]{\shortstack[l]{{\twlrm
9}}}}}
\put(5,125){\makebox(0,0)[lb]{\raisebox{0pt}[0pt][0pt]{\shortstack[l]{{\twlrm
1}}}}}
\put(25,125){\makebox(0,0)[lb]{\raisebox{0pt}[0pt][0pt]{\shortstack[l]{{\twlrm
1}}}}}
\put(45,125){\makebox(0,0)[lb]{\raisebox{0pt}[0pt][0pt]{\shortstack[l]{{\twlrm
1}}}}}
\put(65,125){\makebox(0,0)[lb]{\raisebox{0pt}[0pt][0pt]{\shortstack[l]{{\twlrm
1}}}}}
\put(5,105){\makebox(0,0)[lb]{\raisebox{0pt}[0pt][0pt]{\shortstack[l]{{\twlrm
1}}}}}
\put(25,105){\makebox(0,0)[lb]{\raisebox{0pt}[0pt][0pt]{\shortstack[l]{{\twlrm
1}}}}}
\put(45,105){\makebox(0,0)[lb]{\raisebox{0pt}[0pt][0pt]{\shortstack[l]{{\twlrm
1}}}}}
\put(65,105){\makebox(0,0)[lb]{\raisebox{0pt}[0pt][0pt]{\shortstack[l]{{\twlrm
1}}}}}
\put(85,105){\makebox(0,0)[lb]{\raisebox{0pt}[0pt][0pt]{\shortstack[l]{{\twlrm
2}}}}}
\put(105,105){\makebox(0,0)[lb]{\raisebox{0pt}[0pt][0pt]{\shortstack[l]{{\twlrm
2}}}}}
\put(125,105){\makebox(0,0)[lb]{\raisebox{0pt}[0pt][0pt]{\shortstack[l]{{\twlrm
2}}}}}
\put(305,125){\makebox(0,0)[lb]{\raisebox{0pt}[0pt][0pt]{\shortstack[l]{{\twlrm
1}}}}}
\put(325,125){\makebox(0,0)[lb]{\raisebox{0pt}[0pt][0pt]{\shortstack[l]{{\twlrm
1}}}}}
\put(345,125){\makebox(0,0)[lb]{\raisebox{0pt}[0pt][0pt]{\shortstack[l]{{\twlrm
1}}}}}
\put(365,125){\makebox(0,0)[lb]{\raisebox{0pt}[0pt][0pt]{\shortstack[l]{{\twlrm
1}}}}}
\put(385,125){\makebox(0,0)[lb]{\raisebox{0pt}[0pt][0pt]{\shortstack[l]{{\twlrm
2}}}}}
\put(405,125){\makebox(0,0)[lb]{\raisebox{0pt}[0pt][0pt]{\shortstack[l]{{\twlrm
2}}}}}
\put(425,125){\makebox(0,0)[lb]{\raisebox{0pt}[0pt][0pt]{\shortstack[l]{{\twlrm
2}}}}}
\put(305,105){\makebox(0,0)[lb]{\raisebox{0pt}[0pt][0pt]{\shortstack[l]{{\twlrm
1}}}}}
\put(325,105){\makebox(0,0)[lb]{\raisebox{0pt}[0pt][0pt]{\shortstack[l]{{\twlrm
1}}}}}
\put(345,105){\makebox(0,0)[lb]{\raisebox{0pt}[0pt][0pt]{\shortstack[l]{{\twlrm
1}}}}}
\put(365,105){\makebox(0,0)[lb]{\raisebox{0pt}[0pt][0pt]{\shortstack[l]{{\twlrm
1}}}}}
\put(645,125){\makebox(0,0)[lb]{\raisebox{0pt}[0pt][0pt]{\shortstack[l]{{\twlrm
6}}}}}
\put(665,125){\makebox(0,0)[lb]{\raisebox{0pt}[0pt][0pt]{\shortstack[l]{{\twlrm
6}}}}}
\put(685,125){\makebox(0,0)[lb]{\raisebox{0pt}[0pt][0pt]{\shortstack[l]{{\twlrm
6}}}}}
\put(705,125){\makebox(0,0)[lb]{\raisebox{0pt}[0pt][0pt]{\shortstack[l]{{\twlrm
8}}}}}
\put(725,125){\makebox(0,0)[lb]{\raisebox{0pt}[0pt][0pt]{\shortstack[l]{{\twlrm
8}}}}}
\put(745,125){\makebox(0,0)[lb]{\raisebox{0pt}[0pt][0pt]{\shortstack[l]{{\twlrm
9}}}}}
\put(645,105){\makebox(0,0)[lb]{\raisebox{0pt}[0pt][0pt]{\shortstack[l]{{\twlrm
7}}}}}
\put(665,105){\makebox(0,0)[lb]{\raisebox{0pt}[0pt][0pt]{\shortstack[l]{{\twlrm
7}}}}}
\put(685,105){\makebox(0,0)[lb]{\raisebox{0pt}[0pt][0pt]{\shortstack[l]{{\twlrm
9}}}}}
\put(5,25){\makebox(0,0)[lb]{\raisebox{0pt}[0pt][0pt]{\shortstack[l]{{\twlrm
1}}}}}
\put(25,25){\makebox(0,0)[lb]{\raisebox{0pt}[0pt][0pt]{\shortstack[l]{{\twlrm
1}}}}}
\put(45,25){\makebox(0,0)[lb]{\raisebox{0pt}[0pt][0pt]{\shortstack[l]{{\twlrm
1}}}}}
\put(65,25){\makebox(0,0)[lb]{\raisebox{0pt}[0pt][0pt]{\shortstack[l]{{\twlrm
2}}}}}
\put(5,5){\makebox(0,0)[lb]{\raisebox{0pt}[0pt][0pt]{\shortstack[l]{{\twlrm
1}}}}}
\put(25,5){\makebox(0,0)[lb]{\raisebox{0pt}[0pt][0pt]{\shortstack[l]{{\twlrm
1}}}}}
\put(45,5){\makebox(0,0)[lb]{\raisebox{0pt}[0pt][0pt]{\shortstack[l]{{\twlrm
1}}}}}
\put(65,5){\makebox(0,0)[lb]{\raisebox{0pt}[0pt][0pt]{\shortstack[l]{{\twlrm
1}}}}}
\put(305,25){\makebox(0,0)[lb]{\raisebox{0pt}[0pt][0pt]{\shortstack[l]{{\twlrm
1}}}}}
\put(325,25){\makebox(0,0)[lb]{\raisebox{0pt}[0pt][0pt]{\shortstack[l]{{\twlrm
1}}}}}
\put(345,25){\makebox(0,0)[lb]{\raisebox{0pt}[0pt][0pt]{\shortstack[l]{{\twlrm
1}}}}}
\put(365,25){\makebox(0,0)[lb]{\raisebox{0pt}[0pt][0pt]{\shortstack[l]{{\twlrm
2}}}}}
\put(305,5){\makebox(0,0)[lb]{\raisebox{0pt}[0pt][0pt]{\shortstack[l]{{\twlrm
1}}}}}
\put(325,5){\makebox(0,0)[lb]{\raisebox{0pt}[0pt][0pt]{\shortstack[l]{{\twlrm
1}}}}}
\put(345,5){\makebox(0,0)[lb]{\raisebox{0pt}[0pt][0pt]{\shortstack[l]{{\twlrm
1}}}}}
\put(365,5){\makebox(0,0)[lb]{\raisebox{0pt}[0pt][0pt]{\shortstack[l]{{\twlrm
1}}}}}
\put(565,25){\makebox(0,0)[lb]{\raisebox{0pt}[0pt][0pt]{\shortstack[l]{{\twlrm
1}}}}}
\put(585,25){\makebox(0,0)[lb]{\raisebox{0pt}[0pt][0pt]{\shortstack[l]{{\twlrm
1}}}}}
\put(605,25){\makebox(0,0)[lb]{\raisebox{0pt}[0pt][0pt]{\shortstack[l]{{\twlrm
2}}}}}
\put(625,25){\makebox(0,0)[lb]{\raisebox{0pt}[0pt][0pt]{\shortstack[l]{{\twlrm
3}}}}}
\put(645,25){\makebox(0,0)[lb]{\raisebox{0pt}[0pt][0pt]{\shortstack[l]{{\twlrm
6}}}}}
\put(665,25){\makebox(0,0)[lb]{\raisebox{0pt}[0pt][0pt]{\shortstack[l]{{\twlrm
6}}}}}
\put(685,25){\makebox(0,0)[lb]{\raisebox{0pt}[0pt][0pt]{\shortstack[l]{{\twlrm
6}}}}}
\put(705,25){\makebox(0,0)[lb]{\raisebox{0pt}[0pt][0pt]{\shortstack[l]{{\twlrm
8}}}}}
\put(725,25){\makebox(0,0)[lb]{\raisebox{0pt}[0pt][0pt]{\shortstack[l]{{\twlrm
8}}}}}
\put(745,25){\makebox(0,0)[lb]{\raisebox{0pt}[0pt][0pt]{\shortstack[l]{{\twlrm
9}}}}}
\put(565,5){\makebox(0,0)[lb]{\raisebox{0pt}[0pt][0pt]{\shortstack[l]{{\twlrm
1}}}}}
\put(585,5){\makebox(0,0)[lb]{\raisebox{0pt}[0pt][0pt]{\shortstack[l]{{\twlrm
1}}}}}
\put(605,5){\makebox(0,0)[lb]{\raisebox{0pt}[0pt][0pt]{\shortstack[l]{{\twlrm
1}}}}}
\put(625,5){\makebox(0,0)[lb]{\raisebox{0pt}[0pt][0pt]{\shortstack[l]{{\twlrm
2}}}}}
\put(645,5){\makebox(0,0)[lb]{\raisebox{0pt}[0pt][0pt]{\shortstack[l]{{\twlrm
7}}}}}
\put(665,5){\makebox(0,0)[lb]{\raisebox{0pt}[0pt][0pt]{\shortstack[l]{{\twlrm
7}}}}}
\end{picture}
\end{center}
\normalsize
}
\end{remark}

%%%%%%%%%%%%%%%%%%%%%%%%%%%%%%%%%%%%%%%%%%%%%%%%%%%%%%%%%%%%%%%%%%%%%%%
%%%%%%%%%%%%%%%%%%%%%%%%%%%%%%%%%%%%%%%%%%%%%%%%%%%%%%%%%%%%%%%%%%%%%%%
\subsection{Labelling of cells by ribbon tableaux}

A tabloid $\t$ of shape $\nu=(\nu_1,\ldots,\nu_k)$ can be identified
with a $k$-tuple $(w_1,\ldots,w_k)$ of words, $w_i$ being a row tableau of
lenght $\nu_i$. The Stanton-White correspondence $\psi$ assciates to such
a $k$-tuple of tableaux a $k$-ribbon tableau $T=\psi(\t)$. Thus, the cells
of a unipotent variety $\F^u_\mu$ (where $u$ is of type $\nu$)
are  labelled by $k$-ribbon tableaux of a special kind.  The following
theorem, which implies the stable case of the conjectures, shows that
this labelling is natural from a geometrical point of view.

\begin{theorem}\label{e2cs}
The Stanton-White correspondence $\psi$ sends a tabloid $\t\in L(\nu,\mu)$
onto a ribbon tableau $T=\psi(\t)$ whose cospin is equal to the dimension
of the cell $c_\t$ of $\F^u_\mu$ labelled by $\t$, when one uses the modified
indexation for which the dimension of $c_\t$ is $e(\t)$ (see Section
\ref{HLUV}).
That is,
$$
\cs(\psi(\t))=e(\t) \ .
$$
\end{theorem}

The proof, which is just a direct verification, will not be given here.

At this point, it is useful to observe, following \cite{Te}, that
the $e$-statistic can be given a nonrecursive definition, as a kind
of inversion number. Let $\t=(w_1,\ldots,w_k)$ be a tabloid, identified
with a $k$-tuple of row tableaux.  Let $y$ be the $r$-th letter of $w_i$
and $x$ be the $r$-th letter of $w_j$, and suppose that $x<y$. Then,
the pair $(y,x)$ is said to be an $e$-inversion if either

\medskip
(a) $i<j$

\medskip\noindent or

\medskip
(b) $i>j$ and there s on the right of $x$ in $w_j$ a letter $u<y$

\medskip
Then $e(\t)$ is equal to the number of inversions $(y,x)$ in $\t$.

\begin{example}{\rm Let $\t\in L((2,3,2,1),(2,3,1,1,1))$ be the
following tabloid (the number under a letter $y$ is the number of
$e$-inversions of the form $(y,x)$):
$$
\t \quad = \quad
\left(
\matrix{
\young{2&3\cr} & , &\young{1&1&2\cr}&,&\young{4&5\cr}&,&\young{2\cr} \cr
\matrix{1&1}   &   &\matrix{0&0&0\cr} &&\matrix{3&1\cr}&&\matrix{1\cr} \cr}
\right)
$$
so that $e(\t)=7$. Its image under the SW-correspondence is the $4$-ribbon
tableau

\begin{center}
\setlength{\unitlength}{0.01in}
\begingroup\makeatletter\ifx\SetFigFont\undefined
% extract first six characters in \fmtname
\def\x#1#2#3#4#5#6#7\relax{\def\x{#1#2#3#4#5#6}}%
\expandafter\x\fmtname xxxxxx\relax \def\y{splain}%
\ifx\x\y   % LaTeX or SliTeX?
\gdef\SetFigFont#1#2#3{%
  \ifnum #1<17\tiny\else \ifnum #1<20\small\else
  \ifnum #1<24\normalsize\else \ifnum #1<29\large\else
  \ifnum #1<34\Large\else \ifnum #1<41\LARGE\else
     \huge\fi\fi\fi\fi\fi\fi
  \csname #3\endcsname}%
\else
\gdef\SetFigFont#1#2#3{\begingroup
  \count@#1\relax \ifnum 25<\count@\count@25\fi
  \def\x{\endgroup\@setsize\SetFigFont{#2pt}}%
  \expandafter\x
    \csname \romannumeral\the\count@ pt\expandafter\endcsname
    \csname @\romannumeral\the\count@ pt\endcsname
  \csname #3\endcsname}%
\fi
\fi\endgroup
\begin{picture}(200,95)(0,-10)
\path(0,80)(0,0)
\path(0,0)(200,0)(200,20)
	(200,20)(160,20)(160,40)
	(140,40)(140,80)(0,80)
\path(0,60)(20,60)(20,20)
	(40,20)(40,0)
\path(40,20)(40,80)
\path(120,80)(120,0)
\path(40,60)(120,60)
\path(40,20)(120,20)
\path(60,60)(60,40)(100,40)(100,20)
\path(120,20)(160,20)
\put(5,5){\makebox(0,0)[lb]{\smash{{{\SetFigFont{12}{14.4}{rm}1}}}}}
\put(25,25){\makebox(0,0)[lb]{\smash{{{\SetFigFont{12}{14.4}{rm}2}}}}}
\put(45,45){\makebox(0,0)[lb]{\smash{{{\SetFigFont{12}{14.4}{rm}2}}}}}
\put(65,65){\makebox(0,0)[lb]{\smash{{{\SetFigFont{12}{14.4}{rm}4}}}}}
\put(85,5){\makebox(0,0)[lb]{\smash{{{\SetFigFont{12}{14.4}{rm}1}}}}}
\put(105,25){\makebox(0,0)[lb]{\smash{{{\SetFigFont{12}{14.4}{rm}3}}}}}
\put(125,45){\makebox(0,0)[lb]{\smash{{{\SetFigFont{12}{14.4}{rm}5}}}}}
\put(165,5){\makebox(0,0)[lb]{\smash{{{\SetFigFont{12}{14.4}{rm}2}}}}}
\end{picture}
\end{center}
whose cospin is equal to $7$.
}
\end{example}

\subsection{Atoms}

The $H$-functions indexed by columns can also be completely described
in terms of Hall-Littlewood functions.

\begin{proposition}\label{HCOL}
Let $n=sk+r$ with $0\le r <k$, and set $\lambda=((s+1)^r,s^{k-r})$.
Then,
$$
H^{(k)}_{(1^n)} = \omega \left( \tilde{Q}'_\lambda \right)
$$
where $\omega$ is the involution $s_\lambda\mapsto s_{\lambda'}$.
\end{proposition}

The $k$-quotient
of $(k^n)$ is $(1^s,\ldots,1^s,1^{s+1},\ldots,1^{s+1})$,
where the partition $(1^s)$ is repeated $k-r$ times. Thus,
a $k$-ribbon tableau is mapped by the Stanton-White correspondence
to a $k$-tuple of columns, which can be interpreted as a tabloid,
and the result follows again from Shimomura's decomposition.

\medskip
The partitions arising in Proposition \ref{HCOL} have the property that,
if $\le$ denotes the natural order on partitions,
$$
\alpha \le \lambda \qquad \Longleftrightarrow \qquad \ell(\alpha)\le\ell(\mu) \
{}.
$$
There are canonical injections
$$
\iota_{\alpha\beta} : \tab(\,\cdot\,,\alpha)
\longrightarrow  \tab(\,\cdot\, ,\beta)
$$
when $\alpha\le\beta$ (\cf \cite{LS3,La}). The {\it atom} ${\cal A}(\mu)$
is defined as the set of tableaux in $\tab(\,\cdot\, ,\mu)$ which
are not in the image of any $\iota_{\alpha\mu}$. Define the
symmetric functions ({\it cocharge atoms})
\begin{equation}
\tA_\mu(X;q) =
\sum_\lambda \left(
\sum_{{\bf t}\in {\cal A}(\mu)\cap\tab(\lambda,\mu)} q^{\cc ({\bf t})}
\right)\ s_\lambda(X) \ .
\end{equation}
Proposition \ref{HCOL} can then be rephrased as
\begin{equation}
H^{(k)}_{(1^n)} =
\omega\left(\sum_{\ell(\mu)\le k} \tA_\mu \right) \ .
\end{equation}

It seems that the difference beween the stable $H$-functions
and the immediately lower level can also be described in terms
of atoms. For $\ell(\lambda)=r$, set
$$\tilde D_\lambda(q) = \tilde H^{(r)}_\lambda -\tilde H^{(r-1)}_\lambda \ .$$
These functions seem to be sums of cocharge atoms over
certain intervals in the lattice of partitions.

\begin{conjecture} For any partition $\lambda$, there exists
a partition $f(\lambda)$ such that
$$
\tilde D_\lambda = \sum_{\mu\le f(\lambda)} \tA_\mu \ .
$$
\end{conjecture}

\begin{example}{\rm In weight $6$, the partition $f(\lambda)$
is given by the following table:

\bigskip
\small
\begin{tabular}{|c||c|c|c|c|c|c|c|c|c|c|}\hline
$\lambda$ &
(51)&(42)&(411)&(33)&(321)&(3111)&(222)&(2211)&(21111)&(111111)\\
\hline
$f(\lambda)$ &
(6)&(51)&(51)&(42)&(42)&(411)&(321)&(321)&(3111)&(21111)\\
\hline
\end{tabular}
\normalsize
}
\end{example}

\end{document}